\documentclass[11pt]{article}
\usepackage[latin1]{inputenc}
\usepackage[T1]{fontenc}
\usepackage{amsmath,amssymb}
\usepackage{mathrsfs}
\usepackage[all]{xy}
\usepackage{amsfonts}
\usepackage{amssymb}
\usepackage{lambda}
\usepackage{cite}
\usepackage[oldstyle]{kpfonts}

\def\tab{\hspace*{6mm}}

\def\sgmax{\displaystyle \mathop{\vartriangleleft}_{\hbox{\rm \tiny max}}}

\newcommand{\fin}{
\vskip 2mm
\noindent
$\Box$}

\newcounter{numero}
\newcommand{\num}
{\refstepcounter{numero}
\noindent \thenumero.---}

\begin{document}

\begin{center}
{\LARGE A propos d'une version faible
\vskip 2mm
\noindent
du probl\`eme inverse de Galois}
\vskip 5mm
\noindent
{\large Bruno Deschamps et Fran\c cois Legrand}
\end{center}
\vskip 5mm
\centerline{{\bf R\'esum\'e}}
\vskip 2mm
\noindent
\tab {\small Cet article aborde le Probl\`eme Inverse de Galois Faible qui, pour un corps $k$ donn\'e, \'enonce que, pour tout groupe fini $G$, il existe une extension finie et s\'eparable $L/k$ telle que ${\rm{Aut}}(L/k)=G$. Un de ses objectifs est de montrer comment l'on peut g\'en\'eriquement construire des familles de corps satisfaisant \`a ce probl\`eme, mais pas au traditionnel Probl\`eme Inverse de Galois. C'est par exemple le cas pour les corps $\mathbb{Q}^{\hbox{\scriptsize \rm r\'es}}$, $\mathbb{Q}^{\hbox{\scriptsize \rm tr}}$, $\mathbb{Q}^{\hbox{\scriptsize \rm pyth}}$ ou encore pour les pro-$p$-extensions maximales de $\mathbb{Q}$. Nous d\'emontrons par ailleurs que, pour tout groupe fini non trivial $G$, il existe de nombreux corps v\'erifiant le Probl\`eme Inverse de Galois Faible, mais sur lesquels $G$ ne se r\'ealise pas comme groupe de Galois. Comme autre application, nous montrons que la forme r\'eguli\`ere du Probl\`eme Inverse de Galois Faible admet une r\'eponse positive sur n'importe quel corps.
\vskip 5mm
\centerline{{\bf Abstract}}
\vskip 2mm
\noindent
\tab {\small This paper deals with the Weak Inverse Galois Problem which, for a given field $k$, states that, for every finite group $G$, there exists a finite separable extension $L/k$ such that ${\rm{Aut}}(L/k)=G$. One of its goals is to explain how one can generically produce families of fields which fulfill this problem, but which do not fulfill the usual Inverse Galois Problem. We show that this holds for, e.g., the fields $\mathbb{Q}^{\hbox{\scriptsize \rm sol}}$, $\mathbb{Q}^{\hbox{\scriptsize \rm tr}}$, $\mathbb{Q}^{\hbox{\scriptsize \rm pyth}}$, and for the maximal pro-$p$-extensions of $\mathbb{Q}$. Moreover, we show that, for every finite non-trivial group $G$, there exists many fields fulfilling the Weak Inverse Galois Problem, but over which $G$ does not occur as a Galois group. As a further application, we show that every field fulfills the regular version of the Weak Inverse Galois Problem.}
\noindent
\vskip 10mm
\noindent
{\large \bf 1.--- Introduction.}
\vskip 2mm
\noindent
\tab Le probl\`eme inverse de la th\'eorie de Galois sur un corps $k$ ($\hbox{\rm PIG}_{/k}$ en abr\'eg\'e) consiste \`a savoir si tous les groupes finis apparaissent comme groupes de Galois sur le corps $k$ ou non. Le probl\`eme originel de Hilbert-Noether est le cas $k={\mathbb Q}$ et reste \`a ce jour une question toujours ouverte. L'approche moderne pour tenter de r\'esoudre le $\hbox{\rm PIG}_{/k}$ consiste \`a introduire une ind\'etermin\'ee $T$ et, pour un groupe fini $G$ donn\'e, \`a regarder si l'on peut construire ou non une extension finie galoisienne $E/k(T)$ de groupe de Galois $G$ telle que $E/k$ soit r\'eguli\`ere (c'est-\`a-dire telle que $k$ soit alg\'ebriquement clos dans $E$). Il s'agit du Probl\`eme Inverse de Galois R\'egulier sur $k$ ($\hbox{\rm PIGR}_{/k}$ en abr\'eg\'e), forme g\'eom\'etrique du Probl\`eme Inverse de Galois sur $k$. Nous renvoyons notamment aux livres classiques \cite{Vol96} et \cite{MM99} pour un vaste aper\c cu de ce probl\`eme, ainsi qu'\`a \cite{Zyw14} pour des r\'esultats plus r\'ecents. Il est conjectur\'e (voir par exemple \cite[\S2.1.1]{DD97b}) que le $\hbox{\rm PIGR}_{/k}$ est vrai pour tout corps $k$ (ce qui \'equivaut en fait \`a dire qu'il est vrai sur tout corps premier). Si $k$ est un corps hilbertien\footnote{C'est le cas par exemple si $k$ est un corps de nombres ou le corps des fractions rationnelles en une variable \`a coefficients dans un corps quelconque. Nous renvoyons \`a \cite{FJ08} pour un vaste aper\c cu des corps hilbertiens.}, on voit par sp\'ecialisation que l'on a $\hbox{\rm PIGR}_{/k} \Longrightarrow\hbox{\rm PIG}_{/k}$. Il est donc raisonnable de conjecturer que le $\hbox{\rm PIG}_{/k}$ soit vrai pour tout corps hilbertien $k$, et donc, en particulier, pour $k ={\mathbb Q}$.
\vskip 2mm
\noindent
\tab En 1978, dans \cite{FK78}, E. Fried et J. Koll\'ar ont annonc\'e avoir montr\'e que tout groupe fini apparaissait comme groupe d'automorphismes d'une extension finie (non n\'ecessairement galoisienne) de $\mathbb Q$. Leur preuve comportait cependant une erreur que M. Fried corrigea deux ans plus tard dans \cite{Fri80}. Ce r\'esultat invite naturellement \`a consid\'erer la forme faible du $\hbox{\rm PIG}_{/k}$ suivante : pour tout groupe fini $G$, existe-t-il une extension finie s\'eparable $L/k$ telle que ${\rm{Aut}}(L/k)=G$ ? Dans la suite, nous appellerons cette variante le Probl\`eme Inverse de Galois Faible sur $k$, que nous abr\`egerons en $\hbox{\rm PIGF}_{/k}$. Depuis l'article de M. Fried, plusieurs avanc\'ees significatives sur ce probl\`eme ont \'et\'e effectu\'ees. La plus g\'en\'erale est aussi la plus r\'ecente et est due \`a E. Paran et au second auteur du pr\'esent article qui montrent dans \cite{LP18} que le $\hbox{\rm PIGF}_{/k}$ admet une r\'eponse positive d\`es que $k$ est un corps hilbertien. Ce r\'esultat g\'en\'eralise donc le cas $k={\mathbb Q}$ ainsi que des travaux de T. Takahashi et W.-D. Geyer sur ce probl\`eme\footnote{Nous renvoyons \`a \cite[\S1]{LP18} pour plus de d\'etails et des r\'ef\'erences.}, et vient conforter la conjecture que le $\hbox{\rm PIG}_{/k}$ est vrai pour tout corps hilbertien $k$.
\vskip 2mm
\noindent
\tab \'Etant donn\'e un corps $k$, on peut aussi affaiblir l'\'enonc\'e $\hbox{\rm PIGR}_{/k}$, en demandant seulement de r\'ealiser tout groupe fini $G$ comme le groupe d'automorphismes d'une extension finie s\'eparable $E/k(T)$ telle que $E/k$ soit r\'eguli\`ere. Dans la suite, nous appellerons cette variante le Probl\`eme Inverse de Galois R\'egulier Faible (sur $k$), que nous abr\`egerons en $\hbox{\rm PIGRF}_{/k}$. A notre connaissance, le r\'esultat connu le plus g\'en\'eral sur ce sujet affirme que cet \'enonc\'e poss\`ede une r\'eponse positive pour tout corps $k$ de caract\'eristique nulle, cf. \cite{Fri80}.
\vskip 2mm
\noindent
\tab Le premier \'el\'ement important que nous pr\'esentons dans cet article est que ce r\'esultat est en fait valable pour tout corps. Nous \'etablissons ce fait en nous appuyant sur la preuve originelle de M. Fried et en montrant :
\vskip 2mm
\noindent
{\bf Th\'eor\`eme\num\label{thm intro 1}} {\it Pour tout groupe fini $G$ et tout corps $k$, il existe une extension finie s\'eparable $E/k(T)$ de groupe d'automorphismes $G$ et telle que $E/k$ soit r\'eguli\`ere.}
\vskip 2mm
\noindent
\tab Lorsque le corps de base $k$ est hilbertien, le th\'eor\`eme \ref{thm intro 1} fournit par sp\'ecialisation une r\'eponse positive au $\hbox{\rm PIGF}_{/k}$, ce qui permet de retrouver le r\'esultat principal de \cite{LP18}\footnote{Pour tout corps hilbertien $k$ et tout groupe fini $G$, la m\'ethode utilis\'ee dans \cite{LP18} fournit une extension finie s\'eparable $E/k(T)$ de groupe d'automorphismes $G$ et telle que $E \not \subseteq \overline{k}(T)$. Il s'agit certes d'une conclusion plus faible que celle du th\'eor\`eme \ref{thm intro 1}, mais qui reste bien entendu suffisante pour donner une r\'eponse positive au $\hbox{\rm PIGF}_{/k}$ pour tout corps hilbertien $k$.}. 
\vskip 2mm
\noindent
\tab Du point de vue de la pure th\'eorie des groupes, notre \'etude permet de donner des conditions suffisantes sur l'arithm\'etique d'un corps $k$ pour que le $\hbox{\rm PIGF}_{/k}$ admette une r\'eponse positive. Plus pr\'ecis\'ement, nous montrons que cette conclusion est v\'erifi\'ee s'il existe un groupe fini simple non ab\'elien $G_0$ tel que, pour tous $n \geq 1$ et $N \geq 1$, toute extension\footnote{Dans cet article, nous dirons qu'un groupe $G$ est extension d'un groupe $N$ par un groupe $H$ si l'on a une suite exacte $1 \longrightarrow N \longrightarrow G \longrightarrow H \longrightarrow 1$.} de $G_0^N$ par $A_n$ se r\'ealise comme groupe de Galois sur $k$ (voir lemme \ref{easy} et th\'eor\`eme \ref{thm pigrf}). Ceci nous permet d'exhiber de nombreux corps satisfaisant au PIGF, un premier exemple \'etant donn\'e par le th\'eor\`eme ci-dessous, qui est en fait valable pour tout corps dont le groupe de Galois absolu est librement engendr\'e par une famille infinie d'involutions (voir th\'eor\`eme \ref{1.5 general}).
\vskip 2mm
\noindent
{\bf Th\'eor\`eme\num\label{thm intro 1.5}} {\it Le corps ${\mathbb Q}^{\hbox{\rm \scriptsize tr}}$ des nombres alg\'ebriques totalement r\'eels satisfait au PIGF.}
\vskip 2mm
\noindent
\tab Le r\'esultat principal de \cite{LP18} ne porte que sur les corps hilbertiens et ne permet donc pas de mesurer l'\'eventuelle distance qui pourrait exister entre le $\hbox{\rm PIG}$ et son petit fr\`ere le $\hbox{\rm PIGF}$ (puisque, comme nous l'avons expliqu\'e, ces corps satisfont conjecturalement au PIG). Ce premier exemple est donc int\'eressant dans la mesure o\`u le corps ${\mathbb Q}^{\hbox{\rm \scriptsize tr}}$ ne satisfait pas au $\hbox{\rm PIG}$.
\vskip 2mm
\noindent
\tab Pour aborder le PIGF sur d'autres corps, potentiellement non hilbertiens, nous nous int\'eressons ensuite \`a la notion de cl\^oture d'un corps $k$ relativement \`a une classe ${\mathscr C}$ de groupes finis : il s'agit du corps $k^{\mathscr C}$ qui est, par d\'efinition, le compositum de toutes les extensions finies galoisiennes de $k$ dont le groupe de Galois soit un \'el\'ement de ${\mathscr C}$. Le caract\`ere potentiellement non hilbertien du corps $ k^{\tiny \mathscr C}$ r\'esulte du fait que, sous certaines hypoth\`eses sur la classe ${\mathscr C}$ (cf. th\'eor\`eme \ref{th2brubru}), le corps $k^{\tiny \mathscr C}$ est ${\mathscr C}$-clos, c'est-\`a-dire qu'aucun \'el\'ement non trivial de ${\mathscr C}$ ne se r\'ealise comme groupe de Galois sur $k^{\mathscr C}$. Nous consacrons les sections 3 et 4 de cet article \`a l'\'etude de certaines propri\'et\'es relatives aux classes de groupes finis et \`a celle de l'arithm\'etique des cl\^otures associ\'ees d'un corps. Cette double \'etude, combin\'ee avec celle menant au th\'eor\`eme \ref{thm intro 1}, nous permet alors de montrer le
\vskip 2mm
\noindent
{\bf Th\'eor\`eme\num\label{thm intro 2}} {\it Consid\'erons une classe $\mathscr{C}$ de groupes finis stable par passage au quotient et un corps hilbertien $k$ de caract\'eristique diff\'erente de 2. S'il existe une infinit\'e d'entiers $n \geq 1$ tels que $A_n$ n'appartienne pas \`a $\mathscr{C}$, alors le $\hbox{\rm PIGF}_{/ L}$ admet une r\'eponse positive pour tout corps interm\'ediaire $k \subseteq L \subseteq k^{\tiny \mathscr C}$.}
\vskip 2mm
\noindent
Nous renvoyons au th\'eor\`eme \ref{thm clo1} pour une version plus g\'en\'erale valable aussi en caract\'eristique 2.
\vskip 2mm
\noindent
\tab Nous appliquons enfin le th\'eor\`eme \ref{thm intro 2} \`a plusieurs classes de groupes finis explicites, comme par exemple celle des groupes r\'esolubles (cf. corollaire \ref{coro res}). Une cons\'equence notable de notre \'etude permet alors de pr\'eciser l'\'ecart qui existe entre le Probl\`eme Inverse de Galois et sa version Faible :
\vskip 2mm
\noindent
{\bf Th\'eor\`eme\num\label{thm intro 3}} {\it Pour tout groupe fini non trivial $G$, il existe un corps (non hilbertien) $k$ tel que le $\hbox{\rm PIGF}_{/k}$ admette une r\'eponse positive, mais tel que $G$ ne se r\'ealise pas comme groupe de Galois sur $k$.}
\vskip 2mm
\noindent
On peut en fait obtenir ce m\^eme r\'esultat en prenant \`a la place du seul groupe $G$ une famille finie (et m\^eme infinie sous certaines conditions) de groupes finis non triviaux (voir corollaire \ref{coro P02}).
\vskip 2mm
\noindent
{\bf Remerciements :} Le second auteur b\'en\'eficie des bourses No. 693/13 et 577/15 de l'Israel Science Foundation.
\vskip 10mm
\noindent
{\large \bf 2.--- Le Probl\`eme Inverse de Galois Faible du point de vue de la th\'eorie des groupes.}
\vskip 2mm
\noindent
\tab En toute g\'en\'eralit\'e, si $M/k$ est une extension finie galoisienne de groupe $\Gamma$, la th\'eorie de Galois assure que le groupe d'automorphismes d'une sous-extension quelconque $L/k$ s'identifie au quotient $N_\Gamma(H)/H$, o\`u $H=\hbox{\rm Gal}(M/L)$ et $N_\Gamma(H)$ d\'esigne le normalisateur de $H$ dans $\Gamma$. Cette propri\'et\'e donne un fil conducteur pour aborder, pour un corps $k$ donn\'e, le $\hbox{\rm PIGF}_{/k}$. Ceci nous am\`ene \`a consid\'erer, pour une famille $\{\Gamma_n\}_{n\geq 1}$ de groupes finis, la propri\'et\'e de r\'ealisation suivante :
\vskip 2mm
\noindent
\hbox{\sc (r\'eal)} {\it pour tout groupe fini $G$, il existe un entier $n \geq 1$ et un sous-groupe $H$ de $\Gamma_n$ tel que $N_{\Gamma_n}(H)/H \cong G$.}
\vskip 2mm
\noindent
\tab On voit alors que :
\vskip 2mm
\noindent
{\bf Lemme\num\label{easy}} {\it Pour que le $\hbox{\rm PIGF}_{/k}$ admette une r\'eponse positive, il faut et il suffit qu'il existe une famille $\{\Gamma_n\}_{n\geq 1}$ de groupes finis v\'erifiant la propri\'et\'e de r\'ealisation \hbox{\sc (r\'eal)} et telle que $\Gamma_n$ soit groupe de Galois sur $k$ pour tout $n \geq 1$.}
\vskip2mm
\noindent
\tab Plus g\'en\'eralement, on peut remarquer que, s'il existe une famille $\{\Gamma_n\}_{n\geq 1}$ de groupes finis v\'erifiant la propri\'et\'e de r\'ealisation \hbox{\sc (r\'eal)} et telle que, pour tout $n \geq 1$, il existe une suite $(F_m/k)_{m \geq 1}$ d'extensions finies galoisiennes de groupe de Galois $\Gamma_n$ telle que les corps $F_1, \dots, F_m$ soient lin\'eairement disjoints\footnote{au sens de la d\'efinition donn\'ee \`a la page 35 de \cite{FJ08}.} sur $k$ pour tout $m \geq 2$, alors, pour tout groupe fini $G$, il existe une suite $(E_m/k)_{m \geq 1}$ d'extensions finies s\'eparables de groupe d'automorphismes $G$ telle que $E_1, \dots, E_m$ soient lin\'eairement disjoints sur $k$ pour tout $m \geq 2$.
\vskip2mm
\noindent
\tab L'objectif principal de cette partie est d'\'etablir le r\'esultat suivant, qui fournit de nombreuses situations dans lesquelles la propri\'et\'e de r\'ealisation \hbox{\sc (r\'eal)} est v\'erifi\'ee :
\vskip 2mm
\noindent
{\bf Th\'eor\`eme\num\label{thm pigrf}} {\it Fixons un groupe fini simple non ab\'elien $G_0$ et une suite $(G_n)_{n \geq 1}$ de groupes finis telle que tout groupe fini soit contenu dans au moins un $G_n$. Alors la famille des extensions de $G_0^N$ par $G_n$ (index\'ee par $n \geq 1$ et $N \geq 1$) v\'erifie la propri\'et\'e de r\'ealisation \hbox{\sc (r\'eal)}.}
\vskip 2mm
\noindent
\tab Bien qu'il s'agisse d'un r\'esultat de pure th\'eorie des groupes, nous obtenons en fait celui-ci gr\^ace \`a des arguments de th\'eorie de Galois. Comme premi\`ere application, nous d\'emontrons le th\'eor\`eme \ref{thm intro 1} qui affirme que tout corps satisfait au PIGRF. Ensuite, nous nous int\'eressons aux extensions de groupes dont parle le th\'eor\`eme \ref{thm pigrf}. Nous montrons que, sous une hypoth\`ese ne portant que sur $G_0$, elles sont toujours des produits semi-directs. Ceci nous permet {\it in fine} de prouver le th\'eor\`eme \ref{thm intro 1.5}.
\vskip 5mm
\noindent
{\bf 2.1.--- D\'emonstration du th\'eor\`eme \ref{thm pigrf}.}
\vskip 2mm
\noindent
\tab Pr\'ecisons pour commencer quelques points de terminologie. \'Etant donn\'es un corps $k$ de cl\^oture alg\'ebrique $\overline{k}$ et des ind\'etermin\'ees $T_1,\dots,T_n$, on dira, pour une extension finie galoisienne $E/k(T_1,\dots,T_n)$ donn\'ee, que $E/k$ est {\it r\'eguli\`ere} si $E \cap \overline{k} = k$. Dans le cas $n=1$, si $E \cdot \overline{k}$ d\'esigne le compositum de $E$ et $\overline{k}(T_1)$, un \'el\'ement $t_0$ de $\mathbb{P}^1(\overline{k})$ est un {\it point de branchement} de $E/k(T_1)$ si l'id\'eal premier de $\overline{k}[T_1-t_0]$ engendr\'e par $T_1-t_0$ est ramifi\'e dans l'extension $E \cdot \overline{k}/ \overline{k}(T_1)$ \footnote{Si $t_0= \infty$, $T_1-t_0$ doit \^etre remplac\'e par $1/T_1$.}. Dans la suite, $r$ d\'esignera le nombre de points de branchement de $E/k(T_1)$. Rappelons que $r$ est fini et que l'on a $r=0$ si et seulement si $E \cdot \overline{k}=\overline{k}(T_1)$ (ce qui est \'equivalent \`a $E=k(T_1)$ si $E/k$ est r\'eguli\`ere). Enfin, on dira qu'un groupe fini $G$ est {\it groupe de Galois r\'egulier sur $k$} s'il existe une extension galoisienne $E/k(T_1)$ de groupe $G$ et telle que $E/k$ soit r\'eguli\`ere, et qu'une telle extension $E/k(T_1)$ est une {\it r\'ealisation r\'eguli\`ere de $G$ sur $k$}. 
\vskip2mm
\noindent
\tab La restriction au cas $n=1$ dans la d\'efinition pr\'ec\'edente n'est en fait pas n\'ecessaire, comme le montre le lemme suivant, que nous utiliserons \`a plusieurs reprises dans cet article.
\vskip 2mm
\noindent
{\bf Lemme\num\label{regular}} {\it \'Etant donn\'e un groupe fini $G$, s'il existe des ind\'etermin\'ees $T_1,\dots,T_n$ et une extension galoisienne $E/k(T_1,\dots,T_n)$ de groupe de Galois $G$ et telle que $E/k$ soit r\'eguli\`ere, alors $G$ est groupe de Galois r\'egulier sur $k$.}
\vskip 2mm
\noindent
{\bf Preuve :} La conclusion r\'esulte essentiellement du caract\`ere hilbertien du corps $\overline{k}(T)$, mais on a besoin ici d'une propri\'et\'e un peu plus fine qui assure l'existence de bonnes sp\'ecialisations dans $k(T)$ (et non pas dans $\overline{k}(T)$). Cette propri\'et\'e est classique quand $k$ est infini, cf. par exemple \cite[\S16.2]{FJ08}. Dans le cas o\`u $k$ est fini, nous renvoyons \`a la d\'emonstration de \cite[Lemma 4.2]{DL13} pour plus de d\'etails.
\fin
\vskip 2mm
\noindent
\tab \'Etant donn\'es un groupe fini $G$ et un entier $n \geq 1$ tel que $G$ soit contenu dans $G_n$, nous nous proposons de d\'emontrer le th\'eor\`eme suivant :
\vskip 2mm
\noindent
{\bf{Th\'eor\`eme\num\label{thm pigrf2}}} {\it Fixons un corps $k$ de caract\'eristique $p \geq 0$ tel que l'une des deux conditions suivantes soit satisfaite :
\vskip2mm
\noindent
{\rm{1)}} $G_n$ est groupe de Galois r\'egulier sur $k$ et $G_0$ poss\`ede une r\'ealisation r\'eguli\`ere sur $k$ ayant $r \leq 2$ points de branchement, tous $k$-rationnels,
\vskip2mm
\noindent
{\rm{2)}} $p {\not \vert} 2|G_n|$, $G_n$ est groupe de Galois r\'egulier sur $k$ et $G_0$ poss\`ede une r\'ealisation r\'eguli\`ere sur $k$ ayant $r=3$ points de branchement, tous $k$-rationnels.
\vskip2mm
\noindent
Alors il existe une extension galoisienne $\widetilde{E}/k(T,U)$ telle que $\widetilde{E}/k$ soit r\'eguli\`ere et une sous-extension $E/k(T,U)$ v\'erifiant 
\vskip2mm
\noindent
$\bullet$ ${\rm{Aut}}(E/k(T,U))=G$,
\vskip2mm
\noindent
$\bullet$ ${\rm{Gal}}(\widetilde{E}/k(T,U))$ est une extension $\Gamma$ de $G_0^N$ par $G_n$ pour un certain entier $N \geq 1$.
\vskip 2mm
\noindent
En particulier, si $H$ d\'esigne ${\rm{Gal}}(\widetilde{E}/E)$, on a $N_{\Gamma}(H)/H \cong G$.}
\vskip2mm
\noindent
{\bf Remarque :} Bien qu'il n'y ait aucune hypoth\`ese explicite sur $p$ dans le cas 1), il convient de remarquer que, d'apr\`es le th\'eor\`eme d'existence de Riemann, cette situation ne peut s'appliquer si $p=0$.
\vskip 2mm
\noindent
\tab Avant de d\'etailler la preuve du th\'eor\`eme \ref{thm pigrf2}, nous indiquons comment d\'eduire le th\'eor\`eme \ref{thm pigrf}.
\vskip2mm
\noindent
{\bf Preuve du th\'eor\`eme \ref{thm pigrf} :} Fixons un groupe fini $G$. D'apr\`es la classification des groupes simples finis, $G_0$ est de rang 2. Ainsi, d'apr\`es le th\'eor\`eme d'existence de Riemann, $G_0$ poss\`ede une r\'ealisation r\'eguli\`ere sur $\mathbb{C}$ ayant $r=3$ points de branchement, tous \'etant bien entendu $\mathbb{C}$-rationnels. De plus, pour tout $n \geq 1$ tel que $G \subseteq G_n$, le groupe $G_n$ est groupe de Galois r\'egulier sur $\mathbb{C}$. On obtient alors un groupe $\Gamma$, extension de $G_0^N$ par $G_n$ pour un certain entier $N \geq 1$, et un sous-groupe $H$ de $\Gamma$ tel que $N_{\Gamma}(H)/H \cong G$, par application du cas 2) du th\'eor\`eme \ref{thm pigrf2}\footnote{Notons que le cas 1) du th\'eor\`eme \ref{thm pigrf2} peut \^etre utilis\'e \`a la place du cas 2) (quitte \`a travailler sur $\overline{\mathbb{F}_2}$ au lieu de $\mathbb{C}$) puisque, d'apr\`es \cite{Har84}, tout groupe fini est groupe de Galois r\'egulier sur $\overline{\mathbb{F}_2}$ et, d'apr\`es le th\'eor\`eme de Feit-Thompson et la conjecture d'Abhyankar (d\'emontr\'ee par M. Raynaud et  D. Harbater), tout groupe fini simple non ab\'elien poss\`ede une r\'ealisation r\'eguli\`ere sur $\overline{\mathbb{F}_2}$ n'admettant que $\infty$ comme point de branchement.}.
\fin
\vskip2mm
\noindent
\tab Venons-en maintenant \`a la d\'emonstration du th\'eor\`eme \ref{thm pigrf2}. A partir de maintenant, on se donne deux ind\'etermin\'ees $T$ et $U$, et $\overline{k}$ la cl\^oture alg\'ebrique du corps $k$. On se donne \'egalement une r\'ealisation r\'eguli\`ere $L/k(T)$ de $G_n$ sur $k$. Consid\'erons un \'el\'ement $y(T)$ du sous-corps $L^{G}$ de $L$ tel que $L^{G} = k(T)(y(T))$. On se donne enfin une extension galoisienne $M/k(U)$ de groupe de Galois $G_0$, telle que $M/k$ soit r\'eguli\`ere et telle que $M/k(U)$ poss\`ede
\vskip 2mm
\noindent
$\bullet$ $r \leq 2$ points de branchement, tous $k$-rationnels (si l'on est dans le cas 1)),
\vskip 2mm
\noindent
$\bullet$ $r=3$ points de branchement, tous $k$-rationnels (si l'on est dans le cas 2)).
\vskip 2mm
\noindent
\tab Quitte \`a faire un changement de variable, on peut supposer que l'ensemble des points de branchement de l'extension $M/k(U)$ est
\vskip 2mm
\noindent
$\bullet$ $\{0\}$ si l'on est dans le cas 1) et $r=1$,
\vskip 2mm
\noindent
$\bullet$ $\{0,\infty\}$ si l'on est dans le cas 1) et $r=2$,
\vskip 2mm
\noindent
$\bullet$ $\{0,1,\infty\}$ si l'on est dans le cas 2).
\vskip 2mm
\noindent
\tab Consid\'erons alors un sous-groupe maximal $H$ de $G_0$ et le sous-corps $M^{H}$ de $M$. On se donne un \'el\'ement primitif $x_0(U)$ de $M^{H}/k(U)$, que l'on peut supposer entier sur $k[U]$, et l'on note $P(U,X) \in k[U][X]$ le polyn\^ome minimal de $x_0(U)$ sur $k(U)$. Comme l'extension $M/k$ est r\'eguli\`ere, le polyn\^ome $P(U,X)$ reste irr\'eductible sur $L^{G}(U)$. En particulier, le polyn\^ome tordu $P(U-y(T),X)$ est irr\'eductible sur $L^{G}(U)$. Fixons une racine $x_{y(T)}(U)$ de $P(U-y(T),X)$ et notons $E$ le compositum de $L(U)$ et $L^{G}(U, x_{y(T)}(U))$. Clairement, on a $E=L(U,x_{y(T)}(U))$. 
\vskip 2mm
\noindent
\tab Nous d\'eterminons maintenant le groupe d'automorphismes de l'extension $E/k(T,U)$ :
\vskip 2mm
\noindent
{\bf Lemme\num\label{lemma 2}} {\it On a ${\rm{Aut}}(E/k(T,U))=G$.}
\vskip 2mm
\noindent
{\bf Preuve :} $\bullet$ Commen\c cons par montrer que 
\begin{equation} \label{aut}
{\rm{Aut}}(E/L^{G}(U)) = G.
\end{equation}
Par un argument de r\'egularit\'e d\'ej\`a utilis\'e, le polyn\^ome tordu $P(U-y(T),X)$ est irr\'eductible sur $L(U)$, c'est-\`a-dire les corps $L(U)$ et $L^{G}(U,x_{y(T)}(U))$ sont lin\'eairement disjoints sur $L^{G}(U)$. Comme $L(U)/L^{G}(U)$ est finie galoisienne de groupe $G$, il en est de m\^eme de $E/L^{G}(U,x_{y(T)}(U)).$ Ainsi, pour \'etablir \eqref{aut}, il suffit de voir que tout automorphisme $\sigma$ de $E/L^{G}(U)$ fixe $x_{y(T)}(U)$. Supposons que ce ne soit pas le cas pour un certain $\sigma$. Alors $\sigma(x_{y(T)}(U))$ est une autre racine du polyn\^ome $P(U-y(T),X)$ contenue dans $E$. En particulier, le groupe ${\rm Aut}(E/L(U))$ n'est pas trivial.
\vskip 2mm
\noindent
\tab Consid\'erons maintenant le compositum $M \cdot L$ de $M$ et $L(U)$. Comme l'extension $M/k$ est r\'eguli\`ere, la restriction ${\rm{res}}: {\rm{Gal}}(M \cdot L/L(U)) \longrightarrow {\rm{Gal}}(M/k(U))(=G_0)$ est un isomorphisme de groupes. Notons que le corps $M \cdot L$ est en fait la cl\^oture galoisienne de $L(U,x_0(U))$ sur $L(U)$. Consid\'erons aussi la cl\^oture galoisienne $\widetilde{E}_{y(T)}$ de $E$ sur $L(U)$ et l'automorphisme de corps $\theta \in \hbox{\rm Aut} (L(U)/L)$ d\'efini par $\theta(U) = U-y(T)$. Puisque le polyn\^ome minimal de $x_{y(T)}(U)$ sur $L(U)$ est le tordu par $\theta$ du polyn\^ome minimal de $x_0(U)$ sur $L(U)$, on voit que $\theta$ se rel\`eve en un isomorphisme de corps $\theta:L(U,x_0(U))\longrightarrow L(U,x_{y(T)}(U))$ v\'erifiant  $\theta(x_0(U))=x_{y(T)}(U)$. Un argument classique montre alors que $\theta$ se rel\`eve aux cl\^otures galoisiennes sur $L(U)$ de ces deux corps. Ainsi $\theta$ d\'efinit un isomorphisme de corps $\theta:M \cdot L \longrightarrow \widetilde{E}_{y(T)}$ qui fixe les \'el\'ements de $L$. 
\vskip 2mm
\noindent
\tab Le diagramme suivant r\'ecapitule la situation :
\vskip 2mm
\noindent
{\footnotesize 
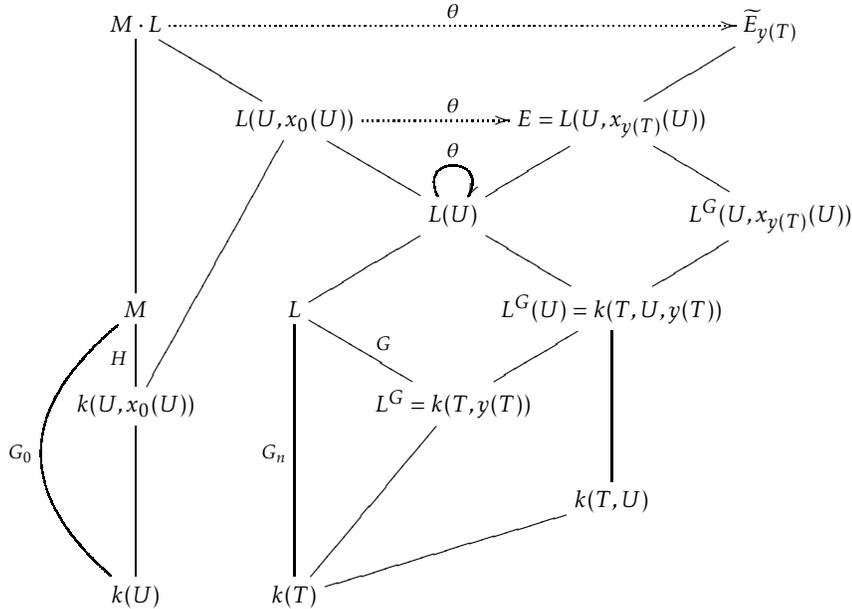
$$\xymatrix @!0 @R=4em @C=5pc {
M \cdot L\ar@{-}[rd]\ar@{.>}[rrrr]^{\theta}& & & & \widetilde{E}_{y(T)}\\
&L(U,x_0(U))\ar@{-}[rd] \ar@{.>}[rr]^{\theta}& & E=L(U, x_{y(T)}(U))\ar@{-}[ur] \ar@{-}[rd] &  \\
&& L(U)\ar@(ul,ur)^{\theta} \ar@{-}[ur] & & L^{G}(U,x_{y(T)}(U)) \\
M\ar@{-}[uuu] \ar@{-} @/_3pc/[ddd]_{G_0}&L \ar@{-}[ur] & & L^{G}(U)=k(T,U,y(T)) \ar@{-}[lu] \ar@{-}[ru]&& \\
k(U,x_0(U))\ar@{-}[u]^{H}\ar@{-}[ruuu]&& L^{G}=k(T,y(T)) \ar@{-}[ru] \ar@{-}[lu]_{G} &    &\\
&& & k(T,U) \ar@{-}[uu] & \\
k(U)\ar@{-}[uu]&k(T) \ar@{-}[rur] \ar@{-}[uuu]^{G_n} \ar@{-}[ruu]& & & \\
}
$$
}
\vskip 2mm
\noindent
\tab L'isomorphisme de corps $\theta$ permet alors de d\'efinir un isomorphisme de groupes $$\psi_{y(T)}:{\rm{Gal}}(\widetilde{E}_{y(T)}/L(U)) \longrightarrow G_0$$ 
en envoyant $\tau \in {\rm{Gal}}(\widetilde{E}_{y(T)}/L(U))$ sur ${\rm{res}}(\theta^{-1} \tau \theta) \in G_0$. De plus, le sous-corps ${\widetilde{E}_{y(T)}}^{\psi_{y(T)}^{-1}(H)}$ de $\widetilde{E}_{y(T)}$ fix\'e par $\psi_{y(T)}^{-1}(H)$ est \'egal \`a $E$. Comme $G_0$ est simple et $H$ est un sous-groupe maximal de $G_0$, on a $N_{G_0}(H)=H$. Ainsi, via $\psi_{y(T)}$, on obtient l'\'egalit\'e $N_{{\rm{Gal}}(\widetilde{E}_{y(T)}/L(U))}({\rm{Gal}}(\widetilde{E}_{y(T)}/E))={\rm{Gal}}(\widetilde{E}_{y(T)}/E)$. Par cons\'equent, ${\rm{Aut}}(E/L(U))$ est trivial, ce qui est impossible en vertu de ce qui pr\'ec\`ede.
\vskip 2mm
\noindent
$\bullet$ Pour \'etablir le lemme, il suffit donc de montrer que ${\rm{Aut}}(E/L^{G}(U))={\rm{Aut}}(E/k(T,U)).$ Il est clair que le groupe de gauche est un sous-groupe de celui de droite. Pour la r\'eciproque, on se donne $\sigma \in {\rm{Aut}}(E/k(T,U))$ et l'on suppose $\sigma \not \in {\rm{Aut}}(E/L^{G}(U))$. Alors $\sigma(y(T)) \not=y(T)$ et $\sigma(x_{y(T)}(U))$ est une racine du polyn\^ome $P(U-\sigma(y(T)),X)$. Comme pr\'ec\'edemment, notons $\widetilde{E}_{y(T)}$ (resp. $\widetilde{E}_{\sigma(y(T))}$) le corps de d\'ecomposition sur $L(U)$ du polyn\^ome $P(U-y(T),X)$ (resp. du polyn\^ome $P(U-\sigma(y(T)),X)$). Par construction, l'ensemble des points de branchement de $\widetilde{E}_{y(T)}/L(U)$ (resp. de $\widetilde{E}_{\sigma(y(T))}/L(U)$) est
\vskip 2mm
\noindent
\tab $\bullet$ $\{y(T)\}$ (resp. $\{\sigma(y(T))\}$) si l'on est dans le cas 1) et $r=1$,
\vskip 2mm
\noindent
\tab $\bullet$ $\{y(T), \infty\}$ (resp. $\{\sigma(y(T)), \infty\}$) si l'on est dans le cas 1) et $r=2$,
\vskip 2mm
\noindent
\tab $\bullet$ $\{y(T), 1+y(T), \infty\}$ (resp. $\{\sigma(y(T)), 1+ \sigma(y(T)), \infty\}$) si l'on est dans le cas 2).
\vskip 2mm
\noindent
Dans chaque cas, on v\'erifie ais\'ement que les extensions $\widetilde{E}_{y(T)}/L(U)$ et $\widetilde{E}_{\sigma(y(T))}/L(U)$ ont des ensembles de points de branchement diff\'erents. En particulier, les corps $\widetilde{E}_{y(T)}$ et $\widetilde{E}_{\sigma(y(T))}$ sont distincts, ce qui entra\^ine $L(U, \sigma(x_{y(T)}(U))) \not=L(U,x_{y(T)}(U))(=E)$. Comme $\sigma$ est dans ${\rm{Aut}}(E/k(T,U))$, on obtient que $L(U, \sigma(x_{y(T)}(U)))$ est strictement contenu dans $E$, ce qui est impossible car ces deux corps sont de degr\'e $|G_0|/|H|$ sur $L(U)$.
\fin
\vskip 2mm
\noindent
\tab Pour tout $\sigma \in {\rm{Gal}}(L(U)/k(T,U))$, notons \`a nouveau $\widetilde{E}_{\sigma(y(T))}$ le corps de d\'ecomposition sur $L(U)$ du polyn\^ome $P(U-\sigma(y(T)),X)$. Comme d\'ej\`a vu, l'extension $\widetilde{E}_{\sigma(y(T))}/L(U)$ est finie galoisienne de groupe de Galois $G_0$. Notons $\widetilde{E}$ le compositum 
$$\bullet_{\sigma} \widetilde{E}_{\sigma(y(T))}$$ 
o\`u $\sigma$ parcourt le groupe ${\rm{Gal}}(L(U)/k(T,U))$. L'extension finie $\widetilde{E}/L(U)$ est galoisienne et, comme $G_0$ est simple, son groupe de Galois est \'egal \`a $G_0^N$ pour un certain entier $N \in \{1, \dots, |G_n|\}$. Remarquons que le corps $E$ construit pr\'ec\'edemmment est contenu dans $\widetilde{E}$ puisque $E$ est contenu dans $\widetilde{E}_{y(T)}$. Par construction, l'extension $\widetilde{E}/k(T,U)$ est galoisienne et son groupe de Galois est une certaine extension de $G_0^N$ par $G_n$. 
\vskip 2mm
\noindent
\tab Pour conclure la d\'emonstration du th\'eor\`eme \ref{thm pigrf2}, il nous reste \`a montrer que l'extension $\widetilde{E}/k$ est r\'eguli\`ere. Par construction, il suffit de d\'emontrer le lemme suivant : 
\vskip 2mm
\noindent
{\bf Lemme\num\label{lem 3}} {\it \'Etant donn\'es un entier $s \geq 2$ et un $s$-uplet $(\sigma_1,\dots, \sigma_s)$ d'\'el\'ements de ${\rm{Gal}}(L(U)/k(T,U))$, supposons que les corps $\widetilde{E}_{\sigma_1(y(T))}, \dots, \widetilde{E}_{\sigma_s(y(T))}$ soient lin\'eairement disjoints sur $L(U)$. Alors les compositums respectifs $$\widetilde{E}_{\sigma_1(y(T))} \cdot \overline{k}, \dots, \widetilde{E}_{\sigma_s(y(T))} \cdot \overline{k}$$ de $\widetilde{E}_{\sigma_1(y(T))}, \dots, \widetilde{E}_{\sigma_s(y(T))}$ et $\overline{k}$ sont lin\'eairement disjoints sur le compositum $L \cdot \overline{k}(U)$ de $L(U)$ et $\overline{k}$.}
\vskip 2mm
\noindent
{\bf Preuve :} Par l'absurde, supposons qu'il existe $q \in \{2,\dots,s\}$ tel que les corps $\widetilde{E}_{\sigma_1(y(T))} \cdot \cdots \cdot \widetilde{E}_{\sigma_{q-1}(y(T))} \cdot \overline{k}$ et $\widetilde{E}_{\sigma_{q}(y(T))} \cdot \overline{k}$ ne soient pas lin\'eairement disjoints sur $L \cdot \overline{k}(U)$. Comme $\widetilde{E}_{\sigma_{q}(y(T))} \cdot \overline{k}/L \cdot \overline{k}(U)$ est de groupe de Galois $G_0$, qui est simple, le corps $\widetilde{E}_{\sigma_{q}(y(T))} \cdot \overline{k}$ est contenu dans le compositum $\widetilde{E}_{\sigma_1(y(T))} \cdot \cdots \cdot \widetilde{E}_{\sigma_{q-1}(y(T))} \cdot \overline{k}$ de $\widetilde{E}_{\sigma_1(y(T))} \cdot \overline{k}, \dots, \widetilde{E}_{\sigma_s(y(T))} \cdot \overline{k}$. En particulier, tout point de branchement de $\widetilde{E}_{\sigma_{q}(y(T))}/L(U)$ est un point de branchement de $\widetilde{E}_{\sigma_1(y(T))} \cdot \cdots \cdot \widetilde{E}_{\sigma_{q-1}(y(T))}/L(U)$. Si l'on est dans le cas 1), cela entra\^ine $\sigma_{q}(y(T))=\sigma_{i}(y(T))$ pour un certain entier $i \in \{1, \dots,q-1\}$, ce qui implique l'\'egalit\'e $\widetilde{E}_{\sigma_{i}(y(T))}=\widetilde{E}_{\sigma_{q}(y(T))}$. En particulier, les corps $\widetilde{E}_{\sigma_1(y(T))} \cdot \cdots \cdot \widetilde{E}_{\sigma_{q-1}(y(T))}$ et $\widetilde{E}_{\sigma_{q}(y(T))}$ ne sont pas lin\'eairement disjoints sur $L(U)$, ce qui est absurde. Si l'on est dans le cas 2), on obtient l'inclusion
$$\{\sigma_{q}(y(T)), 1+\sigma_{q}(y(T))\} \subseteq \bigcup_{i=1}^{q-1} \{\sigma_{i}(y(T)), 1+\sigma_{i}(y(T))\}.$$
Si $\sigma_{q}(y(T))= \sigma_{i}(y(T))$ pour un certain entier $i \in \{1,\dots,q-1\}$, on aboutit comme pr\'ec\'edemment \`a une contradiction. On a donc $\sigma_{q}(y(T)) = 1+ \sigma_i(y(T))$ pour un certain entier $i \in \{1,\dots,q-1\}$. De mani\`ere similaire, on a $1+\sigma_{q}(y(T)) = \sigma_j(y(T))$ pour un certain entier $j \in \{1,\dots,q-1\}$. Les deux \'egalit\'es pr\'ec\'edentes entra\^inent alors $\sigma_j(y(T)) = 2 + \sigma_i(y(T)),$ c'est-\`a-dire $\sigma_j\sigma_i^{-1} (\sigma_i(y(T))) = 2 + \sigma_i(y(T)).$
Par cons\'equent, pour tout entier $l \geq 1$, on a $(\sigma_j\sigma_i^{-1})^l(\sigma_i(y(T))) = 2l + \sigma_i(y(T)).$
Pour $l=|{\rm{Gal}}(L(U)/k(T,U))|=|G_n|$, on obtient $\sigma_i(y(T)) = 2|G_n| + \sigma_i(y(T)).$ Ainsi, $p$ divise $2|G_n|$, ce qui est impossible.
\fin
\vskip 5mm
\noindent
{\bf 2.2.--- Application au Probl\`eme Inverse de Galois R\'egulier Faible.} 
\vskip2mm
\noindent
\tab Nous d\'emontrons maintenant le th\'eor\`eme \ref{thm intro 1}. Fixons pour cela un corps $k$ de caract\'eristique $p \geq 0$ et un groupe fini $G$.
\vskip2mm
\noindent
\tab Supposons tout d'abord $p=2$. Pour $n \geq |G|$, le groupe sym\'etrique $S_n$ contient $G$ et est groupe de Galois r\'egulier sur $k$. De plus, par \cite{AY94b}, le groupe de Mathieu ${\rm{M}}_{23}$ poss\`ede une r\'ealisation r\'eguli\`ere sur $k$ n'admettant que $\infty$ comme point de branchement. Ainsi, en utilisant le lemme \ref{regular} et en appliquant le cas 1) du th\'eor\`eme \ref{thm pigrf2}, on voit qu'il existe une extension finie s\'eparable $E/k(T)$ de groupe d'automorphismes $G$ et telle que $E/k$ soit r\'eguli\`ere.
\vskip2mm
\noindent
\tab Supposons maintenant $p \not=2$. Dans ce cas, nous avons besoin du lemme suivant :
\vskip2mm
\noindent
{\bf Lemme\num\label{An}} {\it{\'Etant donn\'es deux entiers $m \geq 1$ et $n \geq 1$ tels que $m \not \in \{1,2,3,4,6\}$ et $n \not \in \{1,2,3,4,6\}$, ainsi qu'un entier $N \geq 1$, toute extension de $A_m^N$ par $A_n$ est groupe de Galois r\'egulier sur $k$.}}
\vskip2mm
\noindent
{\bf Preuve :} Notons $\Gamma$ une telle extension de $A_m^N$ par $A_n$. Clairement, tout facteur de composition\footnote{On renvoie \`a \cite[\S16.9]{FJ08} pour plus de d\'etails sur la terminologie employ\'ee ici.} de $\Gamma$ est un groupe altern\'e $A_l$ avec $l \not \in \{1,2,3,4,6\}$. Comme $p \not=2$, on peut alors appliquer \cite[Theorem 15]{Bri04} pour affirmer que tout facteur de composition de $\Gamma$ poss\`ede une r\'ealisation GAR sur $k$. Il ne reste alors plus qu'\`a utiliser \cite[\S16.9]{FJ08} et le lemme \ref{regular} pour conclure que $\Gamma$ est groupe de Galois r\'egulier sur $k$.
\fin
\vskip2mm
\noindent
\tab Par le th\'eor\`eme \ref{thm pigrf}, il existe un entier $N \geq 1$, un entier $n \geq 1$ (que l'on peut choisir arbitrairement grand), une extension $\Gamma$ de $A_5^N$ par $A_n$ et un sous-groupe $H$ de $\Gamma$ tel que $N_\Gamma(H)/H \cong G$. Mais, d'apr\`es le lemme \ref{An}, il existe une extension galoisienne $E/k(T)$ de groupe de Galois $\Gamma$ et telle que $E/k$ soit r\'eguli\`ere\footnote{Si $p=0$, cette conclusion peut \^etre obtenue directement gr\^ace au lemme \ref{regular} et au cas 2) du th\'eor\`eme \ref{thm pigrf2}.}. La sous-extension finie s\'eparable $E^H/k(T)$ est alors de groupe d'automorphismes $G$ et l'extension $E^H/k$ est bien entendu r\'eguli\`ere.
\vskip 5mm
\noindent
{\bf 2.3.--- A propos des extensions de groupes apparaissant dans le th\'eor\`eme \ref{thm pigrf}.} 
\vskip2mm
\noindent
\tab La proposition suivante pr\'ecise consid\'erablement les extensions de groupes qui apparaissent dans le th\'eor\`eme \ref{thm pigrf} :
\vskip 2mm
\noindent
{\bf Proposition\num\label{ann3}} {\it Si $G_0$ d\'esigne un groupe fini simple non ab\'elien tel que la suite exacte
$$\xymatrix{1 \ar[r]&\hbox{\rm Int}(G_0) \ar[r]&\hbox{\rm Aut}(G_0)\ar[r]&\hbox{\rm Out}(G_0)\ar@/_1pc/[l]\ar[r]&1}$$
soit scind\'ee, alors toute extension du produit cart\'esien $G_0^N$ ($N\geq 1$) par un groupe fini $G$ est en fait un produit semi-direct $G_0^N \rtimes G$.}
\vskip 2mm
\noindent
\tab La proposition \ref{ann3} d\'ecoule des trois lemmes qui suivent :
\vskip 2mm
\noindent
{\bf Lemme\num\label{ann1}} {\it Si $H$ d\'esigne un groupe de centre trivial et tel que la suite exacte
$$\xymatrix{1 \ar[r]&\hbox{\rm Int}(H) \ar[r]&\hbox{\rm Aut}(H)\ar[r]&\hbox{\rm Out}(H)\ar@/_1pc/[l]\ar[r]&1}$$
soit scind\'ee, alors toute extension de $H$ par un groupe fini $G$ est un produit semi-direct $H \rtimes G$.}
\vskip 2mm
\noindent
{\bf Preuve :} Ce lemme peut \^etre obtenu directement en utilisant la th\'eorie des extensions de groupes d'Eilenberg-MacLane, mais nous allons en donner une preuve \'el\'ementaire (i.e. qui n'utilise pas de cohomologie des groupes).
\vskip 2mm
\noindent
\tab \'Etant donn\'ee une suite exacte $1 \longrightarrow H \longrightarrow \Gamma \longrightarrow G \longrightarrow 1$ avec $G$ fini, notons $s:G\longrightarrow \Gamma$ une section {\it a priori} uniquement \underline{ensembliste}. Consid\'erons le morphisme $I:\Gamma\longrightarrow \hbox{\rm Aut}(H)$ induit par l'\'epimorphisme naturel $\Gamma\longrightarrow  \hbox{\rm Int}(\Gamma)$. Puisque $H$ est de centre trivial, la restriction de $I$ \`a $H$ est un isomorphisme sur $\hbox{\rm Int}(H)$. Par hypoth\`ese, on a une suite exacte scind\'ee
$$\xymatrix{1 \ar[r]&H \ar[r]^-{I}&\hbox{\rm Aut}(H)\ar[r]^{\pi}&\hbox{\rm Out}(H)\ar@/_1pc/[l]_{\sigma}\ar[r]&1}$$
Consid\'erons alors l'application ensembliste
$$\begin{array}{llcl}
\theta:&G&\longrightarrow&\hbox{\rm Aut}(H)\\
&g&\longmapsto&I(s(g)).\\
\end{array}$$
Pour tout $(g_1,g_2)\in G^2$, on a $s(g_1g_2)s(g_2)^{-1}s(g_1)^{-1}\in H$ et donc $\theta(g_1g_2)\theta(g_2)^{-1}\theta(g_1)^{-1}\in \hbox{\rm Int}(H).$ L'application $\overline{\theta}=\pi\circ \theta$ est donc un morphisme $\overline{\theta}: G\longrightarrow \hbox{\rm Out}(H)$. Le diagramme suivant r\'ecapitule la situation :
$$\xymatrix{&H\ar[d]^I\\
&\hbox{\rm Aut}(H)\ar[d]^{\pi}\\
G\ar[r]_{\overline{\theta}}\ar@/_0.5pc/[ur]_{\theta}\ar@/^0.5pc/[ur]^{\sigma\circ \overline{\theta}}&\hbox{\rm Out}(H)\ar@/_1pc/[u]_{\sigma}\\}$$
Puisque $\pi\circ \theta =\pi\circ \sigma \circ \overline{\theta}$, on en d\'eduit que, pour tout $g\in G$, il existe un unique $\lambda(g)\in H$ tel que $\theta (g)=(\sigma \circ \overline{\theta})(g)\circ I(\lambda(g)).$
On consid\`ere alors la section ensembliste $s_0$ de la suite de d\'epart d\'efinie par $s_0(g)=s(g)\lambda(g)^{-1}$ pour tout $g \in G$. Ce qui pr\'ec\`ede montre que, pour tout $g\in G$, on a $\sigma\circ \overline{\theta}(g)=I(s_0(g))$ et donc $I\circ s_0:G\longrightarrow \hbox{\rm Aut}(H)$ est un morphisme. Pour tout $(g_1,g_2)\in G^2$, on a donc $I(s_0(g_1g_2)s_0(g_2)^{-1}s_0(g_1)^{-1})=\hbox{\rm id}_H,$ mais, comme le produit $s_0(g_1g_2)s_0(g_2)^{-1}s_0(g_1)^{-1}$ est un \'el\'ement de $H$ et $H=\hbox{\rm Int}(H)$, on en d\'eduit que $s_0(g_1g_2)s_0(g_2)^{-1}s_0(g_1)^{-1}=1.$ Ceci prouve donc que $s_0$ est un morphisme et que la suite exacte de d\'epart est bien scind\'ee.
\fin
\vskip2mm
\noindent
\tab Dans le cas $N=1$, la proposition \ref{ann3} d\'ecoule de ce lemme. Nous allons maintenant expliquer comment l'obtenir dans le cas g\'en\'eral.
\vskip2mm
\noindent
{\bf Lemme\num\label{ann 1.5}} {\it Fixons un groupe fini simple $G_0$ et un entier $N\geq 1$.
\vskip 2mm
\noindent
{\rm{1)}} Les sous-groupes du produit cart\'esien $G_0^N$ qui sont isomorphes \`a $G_0$ sont exactement ceux de la forme $\left\{(\varphi_1(x),\dots ,\varphi_N(x)) \, \, | \, \,  x\in G_0\right\}$ o\`u $\varphi_1,\dots ,\varphi_N \in \hbox{\rm Aut}(G_0)\cup\{1\}$ avec $(\varphi_1,\dots ,\varphi_N)\ne \{(1,\dots ,1)\}$. 
\vskip 2mm
\noindent
{\rm{2)}} Si $G_0$ est de plus non ab\'elien, alors tout automorphisme $\theta$ de $G_0^N$ op\`ere une permutation des facteurs directs de $G_0^N$, c'est-\`a-dire, si $G_0^N=G_0^{(1)}\times \cdots \times G_0^{(N)}$, alors, pour tout $i \in \{1,\dots, N\}$, il existe $j\in \{1,\dots ,N\}$ tel que $\theta(G_0^{(i)})=G_0^{(j)}$.}
\vskip 2mm
\noindent
{\bf Preuve :} 1) Il est clair que, pour tout $\overline{\varphi}=(\varphi_1,\dots ,\varphi_N)\in (\hbox{\rm Aut}(G_0)\cup\{1\})^N \setminus \{(1,\dots ,1)\}$, le sous-groupe $S_{\overline{\varphi}}=\left\{(\varphi_1(x),\dots ,\varphi_N(x)) \, \, | \, \,  x\in G_0\right\}$ de $G_0^N$ est isomorphe \`a $G_0$. R\'eciproquement, on voit ais\'ement que tout sous-groupe de $G_0^N$ isomorphe \`a $G_0$ est de la forme $S_{\overline{\varphi}}$ avec $\overline{\varphi}=(\varphi_1,\dots ,\varphi_N)\in \hbox{\rm End}(G_0)^N$. Mais, puisque $G_0$ est simple, on a ${\rm{ker}}(\varphi_i)= \{1\}$ ou ${\rm{ker}}(\varphi_i)= G_0$ pour tout $i \in \{1,\dots,N\}$. Dans le premier cas, on a $\varphi_i \in {\rm{Aut}}(G_0)$ alors que, dans le second cas, on a $\varphi_i =1$. Enfin, puisque $G_0$ n'est pas trivial, l'un des $\varphi_i$ est n\'ecessairement \'el\'ement de ${\rm{Aut}}(G_0)$.
\vskip 2mm
\noindent
2) Pour $\overline{\varphi}=(\varphi_1,\dots ,\varphi_N)\in (\hbox{\rm Aut}(G_0)\cup\{1\})^N \setminus \{(1,\dots ,1)$\}, on note $\hbox{\rm Supp}(\overline{\varphi})=\{i \in \{1, \dots, N\} \, \, | \, \, \varphi_i\ne 1\}$($\not= \emptyset$). Puisque $G_0$ est non ab\'elien, si $\overline{\varphi}$ et $\overline{\varphi}'$ sont tels que $\hbox{\rm Supp}(\overline{\varphi})\cap \hbox{\rm Supp}(\overline{\varphi}')\ne \emptyset$, alors $S_{\overline{\varphi}}$ et $S_{\overline{\varphi}'}$ ne commutent pas, c'est-\`a-dire il existe $x\in S_{\overline{\varphi}}$ et $y\in S_{\overline{\varphi}'}$ tels que $xy\ne yx$. On en d\'eduit que, si $S_{\overline{\varphi}_1},\dots ,S_{\overline{\varphi}_N}$ sont des sous-groupes qui commutent deux \`a deux, alors $|\hbox{\rm Supp}(\overline{\varphi}_i)|=1$ pour tout entier $i \in \{1,\dots,N\}$, et donc que les sous-groupes $S_{\overline{\varphi}_1}, \dots, S_{\overline{\varphi}_N}$ sont \'egaux, \`a l'ordre pr\`es, aux facteurs directs $G_0^{(1)}, \dots, G_0^{(N)}$. Si $\theta$ d\'esigne un automorphisme de $G_0^N$, alors les sous-groupes $\theta (G_0^{(1)}),\dots ,\theta (G_0^{(N)})$ sont isomorphes \`a $G_0$ et commutent deux \`a deux. Par le 1) et ce qui pr\'ec\`ede, on d\'eduit que $\theta$ op\`ere bien une permutation des facteurs directs du produit $G_0^N$.
\fin
\vskip 2mm
\noindent
{\bf Lemme\num\label{ann 2}} {\it Pour tout groupe fini simple non ab\'elien $G_0$ et tout entier $N\geq 1$, on a 
$$\hbox{\rm Aut}(G_0^N)=\hbox{\rm Aut}(G_0)^N\rtimes S_N\ \ \hbox{\rm et}\ \ \hbox{\rm Out}(G_0^N)=\hbox{\rm Out}(G_0)^N\rtimes S_N.$$
l'action de $S_N$ sur les groupes $\hbox{\rm Aut}(G_0)^N$ et $\hbox{\rm Out}(G_0)^N$ s'entendant par permutation des facteurs directs. En cons\'equence de quoi, si la suite exacte 
$$\xymatrix{1 \ar[r]&\hbox{\rm Int}(G_0) \ar[r]&\hbox{\rm Aut}(G_0)\ar[r]&\hbox{\rm Out}(G_0)\ar@/_1pc/[l]\ar[r]&1}$$
est scind\'ee, alors la suite exacte 
$$\xymatrix{1 \ar[r]&\hbox{\rm Int}(G_0^N) \ar[r]&\hbox{\rm Aut}(G_0^N)\ar[r]&\hbox{\rm Out}(G_0^N)\ar@/_1pc/[l]\ar[r]&1}$$
l'est aussi.}
\vskip 2mm
\noindent
{\bf Preuve :} Par le 2) du lemme \ref{ann 1.5}, il existe un morphisme $p:\hbox{\rm Aut}(G_0^N)\longrightarrow S_N$ qui, \`a $\theta\in \hbox{\rm Aut}(G_0^N)$, associe la permutation des facteurs directs du produit $G_0^N$ induite par $\theta$. Puisque l'on peut faire op\'erer $S_N$ sur $G_0^N$ en permutant juste les facteurs directs, on voit que ce morphisme est surjectif et que la suite
$$\xymatrix{\hbox{\rm Aut}(G_0^N)\ar[r]^-{p}&S_N\ar@/_1pc/[l]\ar[r]&1}$$
est m\^eme scind\'ee. Si $\theta \in \hbox{\rm ker}(p)$, alors on a $p(G_0^{(i)})=G_0^{(i)}$ pour tout entier $i \in \{1,\dots,N\}$ et donc $\theta\in \hbox{\rm Aut}(G_0)^N\leq \hbox{\rm ker}(p)$. On a donc une suite exacte scind\'ee
$$\xymatrix{1 \ar[r]&\hbox{\rm Aut}(G_0)^N \ar[r]&\hbox{\rm Aut}(G_0^N)\ar[r]&S_N\ar@/_1pc/[l]\ar[r]&1}$$
o\`u l'action de $S_N$ sur $\hbox{\rm Aut}(G_0)^N$ consiste en la permutation des facteurs directs. De plus, comme $\hbox{\rm Int}(G_0^N)=\hbox{\rm Int}(G_0)^N$, on a une suite exacte scind\'ee
$$\xymatrix{1 \ar[r]&\hbox{\rm Out}(G_0)^N \ar[r]&\hbox{\rm Out}(G_0^N)\ar[r]&S_N\ar@/_1pc/[l]\ar[r]&1}$$
o\`u l'action de $S_N$ sur $\hbox{\rm Out}(G_0)^N$ consiste \'egalement en la permutation des facteurs directs. Pour finir, pour chaque entier $i \in \{1,\dots,N\}$, on a par hypoth\`ese une suite exacte scind\'ee
$$\xymatrix@!0 @R=2pc @C=2pc{
1 \ar[rr]&&\hbox{\rm Int}(G_0^{(i)})\ar@{=}[d] \ar[rrr]&&&\hbox{\rm Aut}(G_0^{(i)})\ar@{=}[d]\ar[rrr]&&&\hbox{\rm Out}(G_0^{(i)})\ar@{=}[d]\ar@/_1pc/[lll]\ar[rr]&&1\\
1 \ar[rr]&&\hbox{\rm Int}(G_0)\ar[rrr]&&&\hbox{\rm Aut}(G_0)\ar[rrr]&&&\hbox{\rm Out}(G_0)\ar@/_1pc/[lll]\ar[rr]&&1\\}$$ 
L'action de $S_N$ sur $\hbox{\rm Out}(G_0)^N$ et $\hbox{\rm Aut}(G_0)^N$ \'etant \`a chaque fois la permutation des facteurs, on en d\'eduit (par rel\`evement facteur par facteur) l'existence d'une section \`a la suite exacte
$$\xymatrix{1 \ar[r]&\hbox{\rm Int}(G_0^N) \ar[r]&\hbox{\rm Aut}(G_0^N)\ar[r]&\hbox{\rm Out}(G_0^N)\ar@/_1pc/[l]\ar[r]&1}$$
\fin
\vskip 2mm
\noindent
\tab Combin\'ee avec le th\'eor\`eme \ref{thm pigrf}, la proposition \ref{ann3} montre que la famille de groupes finis $\{A_5^N \rtimes A_n \}$ (index\'ee par $N \geq 1$ et $n \geq 5$) v\'erifie la propri\'et\'e de r\'ealisation \hbox{\sc (r\'eal)} (puisque $\hbox{\rm Int}(A_5)=A_5$, $\hbox{\rm Aut}(A_5)=S_5$ et $\hbox{\rm Out}(A_5)={\mathbb Z}/2\mathbb{Z}$). Par ailleurs, les groupes altern\'es $A_n$ \'etant engendr\'es par des involutions pour $n\geq 5$, on en d\'eduit que :
\vskip 2mm
\noindent
{\bf Th\'eor\`eme\num\label{1.5 general}} {\it Le $\hbox{\rm PIGF}_{/k}$ admet une r\'eponse positive pour tout corps $k$ \`a groupe de Galois absolu librement engendr\'e par une famille infinie d'involutions.}
\vskip 2mm
\noindent
\tab Le th\'eor\`eme \ref{thm intro 1.5} de l'introduction d\'ecoule de ce th\'eor\`eme, compte-tenu du fait que le r\'esultat principal de \cite{FHV93} assure que le groupe de Galois absolu du corps $\mathbb{Q}^{\rm{tr}}$ des nombres alg\'ebriques totalement r\'eels est librement engendr\'e par un ensemble d'involutions hom\'eomorphe \`a l'ensemble triadique de Cantor.
\vskip 10mm
\noindent
{\large \bf 3.--- Quelques propri\'et\'es des classes de groupes finis.}
\vskip 2mm
\noindent
\tab Dans cette section, nous \'etudions quelques propri\'et\'es des classes de groupes finis, en relation avec la th\'eorie des corps.
\vskip 5mm
\noindent
{\bf 3.1.--- Terminologie.}
\vskip 2mm
\noindent
\tab Rappelons tout d'abord (cf. \cite[\S2.1]{RZ10}) qu'une collection non vide ${\mathscr C}$ de groupes est une {\it classe} si ${\mathscr C}$ est stable par isomorphismes, c'est-\`a-dire, si pour tout groupe $G \in {\mathscr C}$ et tout groupe $G'$ isomorphe \`a $G$, on a $G' \in {\mathscr C}$.
\vskip 2mm
\noindent
\tab Dans toute la suite de cet article, on ne consid\`erera que des classes de groupes finis. Pour une telle classe ${\mathscr C}$ donn\'ee, on s'int\'eressera aux quatre propri\'et\'es suivantes :
\vskip 2mm
\noindent
$(\hbox{\rm C}_0)$ $\mathscr C$ est stable par sous-groupes, c'est-\`a-dire, pour tout groupe $G \in {\mathscr C}$ et tout sous-groupe $H$ de $G$, on a $H \in {\mathscr C}$,
\vskip 2mm
\noindent
$(\hbox{\rm C}_1)$ $\mathscr C$ est stable par quotients, c'est-\`a-dire, pour tout groupe $G \in {\mathscr C}$ et tout sous-groupe normal $H$ de $G$, on a $G/H \in {\mathscr C}$,
\vskip 2mm
\noindent
$(\hbox{\rm C}_2)$ $\mathscr C$ est stable par extensions, c'est-\`a-dire, pour tout groupe fini $G$ et tout sous-groupe normal $H$ de $G$ tels que $H$ et $G/H$ soient dans ${\mathscr C}$, on a $G \in {\mathscr C}$,
\vskip 2mm
\noindent
$(\hbox{\rm C}_3)$ $\mathscr C$ est stable par produits fibr\'es surjectifs, c'est-\`a-dire, pour tout produit fibr\'e 
$$\xymatrix @!0 @R=2em @C=4pc {&G_1\ar @{->>}[rd]^{s_1}&\\
G_1\mathop{\times}\limits_{\hbox{\tiny $G_0$}}^{} G_2\ar[ru]\ar[rd]\ar[rr]&&G_0\\
&G_2\ar @{->>}[ru]_{s_2}&\\}
$$
avec $s_1$ et $s_2$ surjectives et tel que $G_1$ et $G_2$ soient dans ${\mathscr C}$, on a $G_1\mathop{\times}\limits_{\hbox{\tiny $G_0$}}^{} G_2\in {\mathscr C}$. 
\vskip 2mm
\noindent
{\bf Remarques :} 1) La propri\'et\'e $(\hbox{\rm C}_3)$ est \'equivalente au fait que, pour tout groupe fini $G$ et tout couple $(H_1,H_2)$ de sous-groupes normaux de $G$ tels que $G/H_1$ et $G/H_2$ soient dans ${\mathscr C}$, on ait $G/(H_1 \cap H_2) \in {\mathscr C}$. En effet, on voit facilement que le groupe quotient $G/(H_1 \cap H_2)$ s'identifie au produit fibr\'e surjectif $G/H_1 \mathop{\times}\limits_{\hbox{\tiny $G/(H_1H_2)$}}^{} G/H_2$.
\vskip 2mm
\noindent
2) La propri\'et\'e $(\hbox{\rm C}_3)$ est v\'erifi\'ee d\`es que les propri\'et\'es $(\hbox{\rm C}_{0,2})$ le sont.
\vskip 2mm
\noindent
3) Les propri\'et\'es $(\hbox{\rm C}_1)$ et $(\hbox{\rm C}_3)$ trouvent chacune une interpr\'etation galoisienne tr\`es claire : les quotients d'un groupe de Galois correspondent aux groupes de Galois des extensions interm\'ediaires et le produit fibr\'e de deux groupes de Galois sur leur intersection correspond au groupe de Galois du compositum des deux extensions associ\'ees.
\vskip 2mm
\noindent
\tab Dans la continuit\'e de la terminologie introduite dans \cite[\S2.1]{RZ10} sur le sujet, nous posons :
\vskip 2mm
\noindent
{\bf D\'efinition\num\label{form-ext}} {\it Nous dirons de la classe ${\mathscr C}$ que c'est 
\vskip 2mm
\noindent
$\bullet$ une "pr\'e-formation" si elle v\'erifie la propri\'et\'e $(\hbox{\rm C}_{1})$,
\vskip 2mm
\noindent
$\bullet$ une "formation" si elle v\'erifie les propri\'et\'es $(\hbox{\rm C}_{1,3})$,
\vskip 2mm
\noindent
$\bullet$ une "formation extensive" si elle v\'erifie les propri\'et\'es $(\hbox{\rm C}_{1,2,3})$,
\vskip 2mm
\noindent
$\bullet$ une "pr\'e-vari\'et\'e" si elle v\'erifie les propri\'et\'es $(\hbox{\rm C}_{0,1})$,
\vskip 2mm
\noindent
$\bullet$ une "vari\'et\'e extensive" si elle v\'erifie les propri\'et\'es $(\hbox{\rm C}_{0,1,2,3})$.}
\vskip 2mm
\noindent
{\bf Exemples :} $\bullet$ Les classes ${\mathscr C}(\hbox{\rm p})$ des $p$-groupes ($p$ premier quelconque), ${\mathscr C}(\hbox{\rm r\'es})$ des groupes finis r\'esolubles, ${\mathscr C}(\hbox{\rm gr})$ de tous les groupes finis et ${\mathscr C}(\hbox{\rm 1})$ compos\'ee uniquement du groupe trivial sont des vari\'et\'es extensives. 
\vskip 2mm
\noindent
$\bullet$ Les classes ${\mathscr C}(\hbox{\rm ab})$ des groupes finis ab\'eliens et ${\mathscr C}(\hbox{\rm nil})$ des groupes finis nilpotents sont \`a la fois des pr\'e-vari\'et\'es et des formations, mais ne sont pas des vari\'et\'es extensives.
\vskip 2mm
\noindent
$\bullet$ La classe ${\mathscr C}(\hbox{\rm cycl})$ des groupes cycliques est une pr\'e-vari\'et\'e qui n'est pas une formation.
\vskip 5mm
\noindent
{\bf 3.2.--- Classe associ\'ee.}
\vskip 2mm
\noindent
\tab Nous g\'en\'eralisons maintenant le passage de la classe des groupes ab\'eliens \`a celle des groupes r\'esolubles.
\vskip 2mm
\noindent
{\bf D\'efinition\num\label{assoc}} {\it On appelle "classe associ\'ee" \`a ${\mathscr C}$ la classe, not\'ee $\widehat{{\mathscr C}}$, des groupes finis $G$ poss\'edant une suite de composition 
$$\{1\} = G_0 \trianglelefteq G_ 1 \trianglelefteq \cdots \trianglelefteq G_{n-1} \trianglelefteq G_n=G$$ telle que les quotients successifs $G_1/G_0,\dots, G_{n}/G_{n-1}$ soient dans ${\mathscr C}$.}
\vskip 2mm
\noindent
{\bf Exemples :} $\bullet$ On a ${\mathscr C}(\hbox{\rm r\'es})=\widehat{{\mathscr C}(\hbox{\rm r\'es})}=\widehat{{\mathscr C}(\hbox{\rm nil})}=\widehat{{\mathscr C}(\hbox{\rm ab})}$.
\vskip 2mm
\noindent
$\bullet$ On a $\widehat{{\mathscr C}(\hbox{\rm p})}={\mathscr C}(\hbox{\rm p})$ pour tout nombre premier $p$.
\vskip 2mm
\noindent
$\bullet$ Puisque tout groupe fini poss\`ede une suite de Jordan-H\"older, on voit que  $\widehat{{\mathscr C}(\hbox{\rm cycl})} = {\mathscr C}(\hbox{\rm r\'es})$ et que, si ${\mathscr C}$ contient la classe ${\mathscr C}(\hbox{\rm simp})$ des groupes finis simples, alors $\widehat{{\mathscr C}}={\mathscr C}(\hbox{\rm gr})$.
\vskip 2mm
\noindent
\tab Nous montrons ci-dessous que certaines des propri\'et\'es pr\'ec\'edemment \'evoqu\'ees se transmettent par passage \`a la classe associ\'ee.
\vskip 2mm
\noindent
{\bf Proposition\num\label{transmission}} {\it {\rm 1)} Si ${\mathscr C}$ v\'erifie $(\hbox{\rm C}_0)$, alors $\widehat{{\mathscr C}}$ v\'erifie $(\hbox{\rm C}_0)$.
\vskip 2mm
\noindent
{\rm 2)} Si ${\mathscr C}$ v\'erifie $(\hbox{\rm C}_1)$, alors $\widehat{{\mathscr C}}$ v\'erifie $(\hbox{\rm C}_1)$.
\vskip 2mm
\noindent
{\rm 3)} La classe $\widehat{{\mathscr C}}$ est la plus petite classe contenant la classe ${\mathscr C}$ et v\'erifiant $(\hbox{\rm C}_2)$.
\vskip 2mm
\noindent
{\rm 4)} Si ${\mathscr C}$ v\'erifie $(\hbox{\rm C}_0)$, alors $\widehat{{\mathscr C}}$ v\'erifie $(\hbox{\rm C}_3)$.}
\vskip 2mm
\noindent
{\bf Preuve :} 1) On se donne $G \in \widehat{{\mathscr C}}$ et un sous-groupe $H$ de $G$. Il existe alors une suite de composition $$\{1\} = G_0 \trianglelefteq G_ 1 \trianglelefteq \cdots \trianglelefteq G_{n-1} \trianglelefteq G_n=G$$ de $G$ telle que les quotients successifs $G_1/G_0,\dots, G_{n}/G_{n-1}$ soient dans ${\mathscr C}$. Consid\'erons la suite de composition
$$\{1\} = G_0 \cap H \trianglelefteq G_1 \cap H \trianglelefteq \cdots \trianglelefteq G_{n-1} \cap H \trianglelefteq G_n \cap H=H$$
de $H$. Pour tout $i \in \{0,\dots,n-1\}$, le quotient $(G_{i+1} \cap H)/(G_i \cap H)$ est isomorphe \`a un sous-groupe de $G_{i+1}/G_i$ et, puisque $G_{i+1}/G_i \in {\mathscr C}$ et ${\mathscr C}$ v\'erifie $(\hbox{\rm C}_0)$, ce sous-groupe est un \'el\'ement de ${\mathscr C}$. Ainsi $H \in \widehat{{\mathscr C}}$.
\vskip 2mm
\noindent
2) On se donne $G \in \widehat{{\mathscr C}}$ et un sous-groupe normal $H$ de $G$. Il existe alors une suite de composition $$\{1\} = G_0 \trianglelefteq G_1 \trianglelefteq \cdots \trianglelefteq G_{n-1} \trianglelefteq G_n=G$$ de $G$ telle que les quotients successifs $G_1/G_0,\dots, G_{n}/G_{n-1}$ soient dans ${\mathscr C}$. Notons $\pi:G \longrightarrow G/H$ la surjection canonique et consid\'erons la suite de composition
$$\{\overline{1}\} = \pi(G_0) \trianglelefteq \pi(G_ 1) \trianglelefteq \cdots \trianglelefteq \pi(G_{n-1}) \trianglelefteq \pi(G_n)=G/H$$
de $G/H$. Pour tout $i \in \{0,\dots,n-1\}$, le quotient $\pi(G_{i+1})/\pi(G_i)$ est isomorphe \`a un quotient de $G_{i+1}/G_i$ et, puisque $G_{i+1}/G_i \in {\mathscr C}$ et ${\mathscr C}$ v\'erifie $(\hbox{\rm C}_1)$, ce quotient est un \'el\'ement de ${\mathscr C}$. Ainsi $G/H \in \widehat{{\mathscr C}}$.
\vskip 2mm
\noindent
3) Montrons tout d'abord que $\widehat{{\mathscr C}}$ v\'erifie $(\hbox{\rm C}_2)$. Pour cela, on se donne un groupe fini $G$ et un sous-groupe normal $H$ de $G$ tels que $H$ et $G/H$ soient dans $\widehat{{\mathscr C}}$. Il existe alors une suite de composition
$$\{1\}=H_0\trianglelefteq H_1 \trianglelefteq \cdots \trianglelefteq H_m=H$$
de $H$ telle que $H_{i+1}/H_i$ soit dans ${\mathscr C}$ pour tout $i \in \{0,\dots ,m-1\}$, et une suite $$H=H_{m}\trianglelefteq H_{m+1} \trianglelefteq \cdots \trianglelefteq H_n=G$$ telle que $H\leq H_i$ pour tout $i \in \{m,\dots ,n\}$ et telle que
$$\{\overline{1}\}=H_m/H\trianglelefteq H_{m+1}/H \trianglelefteq \cdots \trianglelefteq H_n/H=G/H$$
soit une suite de composition de $G/H$ v\'erifiant $(H_{i+1}/H)/(H_i/H)\in {\mathscr C}$ pour tout $i \in \{m,\dots ,n-1\}$. Puisque $(H_{i+1}/H)/(H_i/H)\simeq  H_{i+1}/H_i$ pour tout $i \in \{m,\dots ,n-1\}$, on en d\'eduit que
$$\{1\}=H_0\trianglelefteq \cdots \trianglelefteq H_m=H\trianglelefteq H_{m+1}\trianglelefteq \cdots \trianglelefteq H_n=G$$
est une suite de composition de $G$ telle que $H_{i+1}/H_i \in {\mathscr C}$ pour tout $i \in \{0,\dots ,n-1\}$. Ainsi $G \in \widehat{{\mathscr C}}$.
\vskip 2mm
\noindent
\tab On se donne maintenant une classe ${\mathscr C}_0$ de groupes finis contenant ${\mathscr C}$ et v\'erifiant $(\hbox{\rm C}_2)$. Fixons un groupe fini $G\in \widehat{\mathscr C}$, muni d'une suite de composition
$$\{1\} = G_0 \trianglelefteq G_1 \trianglelefteq \cdots \trianglelefteq G_{n-1} \trianglelefteq G_n=G$$
telle que les quotients successifs $G_1/G_0,\dots, G_{n}/G_{n-1}$ soient dans ${\mathscr C}$. Comme $G_0$ est trivial et $G_1/G_0 \in {\mathscr C}$, on voit que $G_1 \in {\mathscr C} \subseteq {\mathscr C}_0$. De plus, $G_2/G_1 \in {\mathscr C} \subseteq {\mathscr C}_0$ et ${\mathscr C}_0$ v\'erifie $(\hbox{\rm C}_2)$. Ainsi $G_2$ est un \'el\'ement de ${\mathscr C}_0$. Par r\'ecurrence, on voit alors que $G_n=G$ est un \'el\'ement de ${\mathscr C}_0$. Ainsi $\widehat{\mathscr C}\subseteq {\mathscr C}_0$.
\vskip 2mm
\noindent
4) Il s'agit d'une cons\'equence imm\'ediate du 1) et du 3) ci-dessus, et du 2) des remarques du \S3.1.
\fin
\vskip 2mm
\noindent
{\bf Corollaire\num\label{ffr}} {\it Si ${\mathscr C}$ est une pr\'e-vari\'et\'e, alors $\widehat{{\mathscr C}}$ est une vari\'et\'e extensive.}
\vskip 5mm
\noindent
{\bf 3.3.--- Classe duale.}
\vskip 2mm
\noindent
\tab Nous \'etudions maintenant une notion capitale pour notre propos : celle de classe duale.
\vskip 2mm
\noindent
{\bf D\'efinition\num} {\it On appelle "classe duale" de ${\mathscr C}$ la classe, not\'ee ${\mathscr C}^*$, des groupes finis $G$ tels que, pour tout sous-groupe normal $H$ de $G$ tel que $G/H$ soit dans ${\mathscr C}$, on ait $H=G$.}
\vskip 2mm
\noindent
\tab Cette notion est intimement li\'ee \`a la th\'eorie inverse de Galois : si ${\mathscr C}$ d\'esigne la classe des groupes finis qui n'apparaissent pas comme groupes de Galois sur un corps $k$ fix\'e, on voit que, pour qu'un groupe fini $G$ donn\'e soit groupe de Galois sur $k$, il faut n\'ecessairement que $G$ soit \'el\'ement de ${\mathscr C}^*$. L'\'etude de la classe duale va jouer un r\^ole important dans cet article, en particulier l'\'etude des propri\'et\'es de ${\mathscr C}^*$ qui d\'ecoulent de celles de ${\mathscr C}$. Il est d\'ej\`a facile de voir que le passage \`a la classe duale dualise certaines des op\'erations ensemblistes usuelles : si ${\mathscr C}_1$ et ${\mathscr C}_2$ d\'esignent deux classes quelconques de groupes finis, alors
\vskip 2mm
\noindent
a) $({\mathscr C}_1\cup {\mathscr C}_2)^*={\mathscr C}^*_1\cap {\mathscr C}_2^*$,
\vskip 2mm
\noindent
b) ${\mathscr C}_1\subseteq {\mathscr C}_2 \Longrightarrow {\mathscr C}_2^*\subseteq {\mathscr C}_1^*$,
\vskip 2mm
\noindent
c) ${\mathscr C}_1^*\cup {\mathscr C}_2^* \subseteq ({\mathscr C}_1\cap {\mathscr C}_2)^*$ \footnote{L'inclusion r\'eciproque est fausse en g\'en\'eral. En effet, consid\'erons par exemple les classes ${\mathscr C}_1=\{\{1\},{\mathbb Z}/2\mathbb{Z}\}$ et ${\mathscr C}_2=\{\{1\},{\mathbb Z}/3{\mathbb Z}\}$. On a alors ${\mathscr C}_1\cap {\mathscr C}_2={\mathscr C}(\hbox{\rm 1})$ et donc $({\mathscr C}_1\cap {\mathscr C}_2)^*={\mathscr C}(\hbox{\rm gr})$, alors que ${\mathbb Z}/6{\mathbb Z} \notin {\mathscr C}_1^*\cup {\mathscr C}_2^*$.}.
\vskip 2mm
\noindent
\tab De mani\`ere un peu moins triviale, on a :
\vskip 2mm
\noindent
{\bf Proposition\num\label{transmission 2}} {\it {\rm 1)} La classe ${\mathscr C}^*$ v\'erifie $(\hbox{\rm C}_1)$.
\vskip 2mm
\noindent
{\rm 2)} Si ${\mathscr C}$ v\'erifie $(\hbox{\rm C}_1)$, alors ${\mathscr C}^*$ v\'erifie $(\hbox{\rm C}_2)$.
\vskip 2mm
\noindent
{\rm 3)} On a ${\mathscr C}^*=\big({\widehat{\mathscr C}}\big)^{*}$.}
\vskip 2mm
\noindent
{\bf Preuve :} 1) On se donne un groupe fini $G \in {\mathscr C}^*$, un sous-groupe normal $H$ de $G$ et un sous-groupe normal $V$ de $G/H$ tel que $(G/H)/V \in {\mathscr C}$. Clairement, on a $V=H'/H$ pour un certain sous-groupe normal $H'$ de $G$ contenant $H$. Comme $(G/H)/V \in {\mathscr C}$ et $(G/H)/V = (G/H)/(H'/H) \simeq G/H'$, on a $G/H' \in {\mathscr C}$. Le groupe $G$ \'etant dans ${\mathscr C}^*$, on obtient $G=H'$, c'est-\`a-dire $V=G/H$. Ainsi $G/H \in {\mathscr C}^*$.
\vskip 2mm
\noindent
2) On se donne un groupe fini $G$ et un sous-groupe normal $H$ de $G$ tels que $H$ et $G/H$ soient dans ${\mathscr C}^*$. Par l'absurde, supposons que $G$ ne soit pas dans ${\mathscr C}^*$. Il existe alors un sous-groupe normal $H'$ de $G$ tel que $H' \not=G$ et $G/H' \in {\mathscr C}$. S'il existe un sous-groupe normal $H''$ de $G$ tel que $H' \leq H'' \not=G$, on a $G/H'' \simeq (G/H')/(H''/H')$. Comme ${\mathscr C}$ v\'erifie $(\hbox{\rm C}_1)$ et $G/H'$ est dans ${\mathscr C}$, on voit que $G/H''$ est lui aussi dans ${\mathscr C}$. On peut donc supposer que $G/H'$ est simple. Si $H$ est contenu dans $H'$, alors $H'/H$ est un sous-groupe normal de $G/H$ et $(G/H)/(H'/H) \cong G/H' \in {\mathscr C}$. Or $G/H$ est un \'el\'ement de ${\mathscr C}^*$ par hypoth\`ese. On a donc $H'/H=G/H$, c'est-\`a-dire $H'=G$, ce qui est absurde. Si $H$ n'est pas contenu dans $H'$, alors $HH'$ est un sous-groupe normal de $G$ contenant strictement $H'$. Comme $G/H'$ est simple, on a $HH'=G$. Ainsi $H/H \cap H' \simeq HH'/H' = G/H' \in {\mathscr C}$. Or $H$ est un \'el\'ement de ${\mathscr C}^*$ par hypoth\`ese. On a donc $H \cap H' = H$, ce qui est absurde. On en d\'eduit ainsi que $G$ est \'el\'ement de ${\mathscr C}^*$.
\vskip 2mm
\noindent
3) Puisque ${\mathscr C} \subseteq \widehat{{\mathscr C}}$, on a $\big(\widehat{\mathscr C}\big)^*\subseteq {\mathscr C}^*$. R\'eciproquement, on se donne un groupe fini $G$ n'appartenant pas \`a $\big(\widehat{\mathscr C}\big)^*$. Il existe alors un sous-groupe normal $H$ de $G$ tel que $H \not=G$ et $G/H \in \widehat{{\mathscr C}}$. Par d\'efinition de $\widehat{\mathscr C}$, il existe un sous-groupe normal $H'$ de $G/H$ tel que $(G/H)/H' \in {\mathscr C}$, et que l'on peut de plus supposer diff\'erent de $G/H$ puisque $G/H \not= \{1\}$. Clairement, on a $H'=H''/H$ pour un certain sous-groupe normal $H''$ de $G$ contenant $H$. Ainsi $G/H'' \in {\mathscr C}$ et on a $H'' \not=G$. Par cons\'equent, $G \not \in {\mathscr C}^*$.
\fin
\vskip 5mm
\noindent
{\bf 3.4.--- Multidualit\'e.}
\vskip 2mm
\noindent
\tab Nous finissons cette section en nous int\'eressant aux classes duales successives de la classe ${\mathscr C}$. Commen\c cons par l'\'etude de la {\it{classe biduale}} ${\mathscr C}^{**}=({\mathscr C}^{*})^{*}$ de ${\mathscr C}$. Dans ce qui suit, on notera $H\vartriangleleft G$ pour dire que $H$ est un sous-groupe normal strict de $G$ et $H\sgmax G$ pour dire que $H$ est maximal dans l'ensemble des sous-groupes normaux stricts de $G$ (c'est-\`a-dire $G/H$ est simple). 
\vskip 2mm
\noindent
\tab Pour un groupe fini $G$ donn\'e, on a :
$$\begin{array}{lll}
G\in {\mathscr C}^{**}&\Longleftrightarrow &\forall H\vartriangleleft G,\ G/H\notin {\mathscr C}^{*}\\
&\Longleftrightarrow &\forall H\vartriangleleft G,\ \exists \overline{N}\vartriangleleft G/H,\ (G/H)/\overline{N} \in {\mathscr C}\\
&\Longleftrightarrow &\forall H\vartriangleleft G,\ \exists N\vartriangleleft G\ \hbox{\rm tel que}\ H\leq N\ \hbox{\rm et}\  G/N \in {\mathscr C}\\
&\Longleftrightarrow &\forall H\sgmax G,\ G/H \in {\mathscr C}.\\
\end{array}$$
\vskip 2mm
\noindent
\tab De cette caract\'erisation d\'ecoulent les deux remarques suivantes :
\vskip 2mm
\noindent
a) Si ${\mathscr C}$ v\'erifie $(\hbox{\rm C}_1)$, alors ${\mathscr C}$ est contenue dans ${\mathscr C}^{**}$. En g\'en\'eral, l'inclusion ${\mathscr C} \subseteq {\mathscr C}^{**}$ n'est pas vraie, comme le montre l'exemple ${\mathscr C}=\{{\mathbb Z}/4{\mathbb Z}\}$.
\vskip 2mm
\noindent
b) Si ${\mathscr C}({\rm simp}) \subseteq {\mathscr C}$, alors ${\mathscr C}^{**}={\mathscr C}({\rm gr})$. Cette remarque prouve en particulier que l'inclusion r\'eciproque ${\mathscr C}^{**} \subseteq {\mathscr C}$ n'est pas vraie non plus en g\'en\'eral.
\vskip 2mm
\noindent
\tab Pour un groupe fini non trivial $G$ donn\'e, rappelons que le {\it radical de Baer} de $G$ (cf. \cite{Bae64}), que nous noterons $\hbox{\rm Rad}(G)$, est l'intersection de tous les sous-groupes normaux maximaux de $G$ :
$$\hbox{\rm Rad}(G)=\bigcap_{H \sgmax G} H.$$
\vskip 1mm
\noindent
\tab Nous pouvons maintenant caract\'eriser les \'el\'ements non triviaux de ${\mathscr C}^{**}$ de la mani\`ere suivante :
\vskip 2mm
\noindent
{\bf Proposition\num\label{prop bidual}} {\it \'Etant donn\'e un groupe fini non trivial $G$, les deux assertions 
\vskip 2mm
\noindent
{\rm i)} $G \in {\mathscr C}^{**}$,
\vskip 2mm
\noindent
{\rm ii)} $G/\hbox{\rm Rad}(G)$ est isomorphe \`a un produit direct non vide de groupes finis simples appartenant \`a ${\mathscr C}$,
\vskip 2mm
\noindent
sont \'equivalentes.}
\vskip 2mm
\noindent
{\bf Preuve :} Parmi les sous-groupes normaux $H$ de $G$ tels que $H \sgmax G$, l'on consid\`ere une famille finie $\{H_1, \dots, H_n\}$ telle que $\hbox{\rm Rad}(G)=H_1\cap \cdots \cap H_n$ et telle que $H_{i+1}\cap \cdots \cap H_n \not \subseteq H_i$ pour tout $i \in \{1,\dots ,n-1\}$.
\vskip 2mm
\noindent
\tab Comme $G/H_1$ est simple et $K=H_2\cap \cdots \cap H_n$ est un sous-groupe normal de $G$ non contenu dans $H_1$, on a $H_1 K=G$. Ainsi, en appliquant le deuxi\`eme th\'eor\`eme d'isomorphisme, on a
$$|G/(H_1 \cap K)| = \frac{|H_1 K|}{|H_1 \cap K|} = \frac{|H_1||K|}{|H_1 \cap K|^2} = \frac{|H_1K|}{|H_1|}\frac{|H_1K|}{|K|}$$
et donc $|G/(H_1 \cap K)|=|G/H_1 \times G/K|$. Par ailleurs, le morphisme canonique $G \longrightarrow G/H_1 \times G/K$ est de noyau $H_1 \cap K$, ce qui montre finalement que 
$$G/{\hbox{\rm Rad}}(G) = G/(H_1 \cap K) \simeq G/H_1 \times G/K = G/H_1 \times G/(H_2 \cap \cdots \cap H_n).$$ 
Une r\'ecurrence imm\'ediate montre alors que
$G/\hbox{\rm Rad}(G) \simeq G/H_1\times \cdots \times G/H_n.$
\vskip 2mm
\noindent
\tab Venons-en maintenant \`a la d\'emonstration de l'\'equivalence annonc\'ee. Supposons tout d'abord que $G$ soit dans ${\mathscr C}^{**}$. D'apr\`es ce qui pr\'ec\`ede, on a $G/\hbox{\rm Rad}(G) \simeq G/H_1\times \cdots \times G/H_n$, et chaque quotient $G/H_i$ est un groupe simple appartenant \`a ${\mathscr C}$, en vertu de la caract\'erisation des \'el\'ements de la classe biduale vue pr\'ec\'edemment. R\'eciproquement, supposons que $G/{\hbox{\rm Rad}}(G)$ soit isomorphe \`a un produit direct $G_1 \times \cdots \times G_r$ de groupes simples appartenant \`a ${\mathscr C}$. Pour tout $H \sgmax G$, on a alors un \'epimorphisme $G_1 \times \cdots \times G_r \longrightarrow G/H$. Comme les groupes $G/H, G_1, \dots, G_r$ sont simples, cet \'epimorphisme fournit un isomorphisme $G_j \cong G/H$ pour un certain $j \in \{1,\dots,r\}$. Ainsi $G/H\in{\mathscr C}$ et $G$ est donc dans ${\mathscr C}^{**}$.
\fin
\vskip 2mm
\noindent
\tab D\'efinissons maintenant la {\it{$n$-i\`eme classe duale}} ${\mathscr C}^{*[n]}$ de ${\mathscr C}$ par r\'ecurrence sur l'entier $n$, en posant ${\mathscr C}^{*[0]}={\mathscr C}$ et, pour tout $n\geq 0$, ${\mathscr C}^{*[n+1]}=\big({\mathscr C}^{*[n]} \big)^*$.
\vskip 2mm
\noindent
\tab Nous avons alors le r\'esultat de cyclicit\'e des classes multiduales suivant :
\vskip 2mm
\noindent
{\bf Proposition \num\label{multi}} {\it {\rm 1)} On a ${\mathscr C}^* \subseteq {\mathscr C}^{***}$ et, si ${\mathscr C}$ v\'erifie $(\hbox{\rm C}_1)$, alors ${\mathscr C}^{*}={\mathscr C}^{***}$.
\vskip 2mm
\noindent
{\rm 2)} Pour tout entier $n\geq 1$, on a
${\mathscr C}^{*[2n]}={\mathscr C}^{**}$ et ${\mathscr C}^{*[2n+1]}={\mathscr C}^{***}$.}
\vskip 2mm
\noindent
{\bf Preuve :} 1) D'apr\`es le 1) de la proposition \ref{transmission 2}, la classe ${\mathscr C}^*$ v\'erifie $(\hbox{\rm C}_1)$ et on a donc ${\mathscr C}^* \subseteq {\mathscr C}^{***}$ en vertu de la remarque a) ci-dessus. Si ${\mathscr C}$ v\'erifie $(\hbox{\rm C}_1)$, alors on a ${\mathscr C} \subseteq {\mathscr C}^{**}$, et donc ${\mathscr C}^{***} \subseteq {\mathscr C}^{*}$.
\vskip 2mm
\noindent
2) Il suffit d'effectuer une r\'ecurrence sur l'entier $n$ en appliquant le 1) et en remarquant que, d'apr\`es le 1) de la proposition \ref{transmission 2}, la classe ${\mathscr C}^{*[n]}$ v\'erifie $(\hbox{\rm C}_1)$ pour tout $n \geq 1$.
\fin
\vskip 2mm
\noindent
\tab L'\'etude de la multidualit\'e de la classe ${\mathscr C}$ se r\'esume donc \`a la donn\'ee de ${\mathscr C}^*$, ${\mathscr C}^{**}$ et ${\mathscr C}^{***}$ (et, si ${\mathscr C}$ v\'erifie $(\hbox{\rm C}_1)$, seulement \`a la donn\'ee des classes duale et biduale de ${\mathscr C}$). L'exemple de la classe ${\mathscr C}=\{\{1\},{\mathbb Z}/4{\mathbb Z}\}$ montre que ${\mathscr C}$, ${\mathscr C}^*\notin \{{\mathscr C}, {\mathscr C}({\rm 1}),{\mathscr C}({\rm gr})\}$, ${\mathscr C}^{**}={\mathscr C}({\rm 1})$ et ${\mathscr C}^{***}={\mathscr C}({\rm gr})$ peuvent \^etre des classes deux \`a deux distinctes.
\vskip 10mm
\noindent
{\large \bf 4.--- Arithm\'etique des $\mathscr C$-cl\^otures d'un corps.}
\vskip 2mm
\noindent
\tab On s'int\'eresse maintenant aux extensions galoisiennes finies \`a groupe de Galois dans une classe donn\'ee. Dans toute cette section, $k$ d\'esigne un corps de cl\^oture s\'eparable $k^{\hbox{\rm \scriptsize s\'ep}}$ et ${\mathscr C}$ une classe de groupes finis.
\vskip 2mm
\noindent
{\bf D\'efinition\num\label{def}} {\it On appelle "$\mathscr C$-cl\^oture" de $k$ le corps, not\'e $k^{\tiny \mathscr C}$, \'egal au compositum (dans $k^{\hbox{\rm \scriptsize s\'ep}}$) de toutes les extensions galoisiennes finies de $k$ ayant un groupe de Galois \'el\'ement de $\mathscr C$. Le corps $k$ est dit "$\mathscr C$-clos" si $k=k^{\tiny \mathscr C}$, c'est-\`a-dire si $k$ ne poss\`ede aucune extension galoisienne finie non triviale \`a groupe de Galois dans ${\mathscr C}$.}
\vskip 5mm
\noindent
{\bf 4.1.--- Groupes de Galois et ${\mathscr C}$-cl\^otures.}
\vskip 2mm 
\noindent
\tab Si l'on prend pour ${\mathscr C}$ la classe ${\mathscr C}(\hbox{\rm cycl})$ ou la classe ${\mathscr C}(\hbox{\rm ab})$, le corps ${\mathbb Q}^{\tiny \mathscr C}$ est alors la traditionnelle cl\^oture cyclotomique de $\mathbb Q$, souvent not\'ee ${\mathbb Q}^{\hbox{\rm \scriptsize ab}}$. Tous les groupes ab\'eliens finis se r\'ealisant comme groupes de Galois sur $\mathbb Q$, l'on voit que, bien que ${\mathbb Q}^{\hbox{\rm \scriptsize ab}}$ soit la ${\mathscr C}(\hbox{\rm cycl})$-cl\^oture de $\mathbb Q$, l'extension ${\mathbb Q}^{\hbox{\rm \scriptsize ab}}/{\mathbb Q}$ poss\`ede des sous-extensions finies galoisiennes sur $\mathbb{Q}$ \`a groupes de Galois non cycliques. Cette pathologie dispara\^it si l'on regarde le corps ${\mathbb Q}^{\hbox{\rm \scriptsize ab}}$ comme la ${\mathscr C}(\hbox{\rm ab})$-cl\^oture de $\mathbb Q$. Pour autant, on peut facilement construire des extensions finies ab\'eliennes non triviales de ${\mathbb Q}^{\hbox{\rm \scriptsize ab}}$, et ${\mathbb Q}^{\hbox{\rm \scriptsize ab}}$ n'est donc pas ${\mathscr C}(\hbox{\rm ab})$-clos. La ${\mathscr C}(\hbox{\rm r\'es})$-cl\^oture de ${\mathbb Q}$, not\'ee ${\mathbb Q}^{\hbox{\rm \scriptsize r\'es}}$ pour plus de commodit\'e, ne pr\'esente elle aucune de ces pathologies : les sous-extensions galoisiennes finies sur $\mathbb{Q}$ de ${\mathbb Q}^{\hbox{\rm \scriptsize r\'es}}$ ont toutes un groupe de Galois r\'esoluble et aucune extension galoisienne finie non triviale de ${\mathbb Q}^{\hbox{\rm \scriptsize r\'es}}$ ne poss\`ede un groupe de Galois r\'esoluble. Le th\'eor\`eme qui suit vise \`a d\'eterminer quelles propri\'et\'es il faut demander \`a la classe ${\mathscr C}$ pour que la ${\mathscr C}$-cl\^oture associ\'ee ait ces "bonnes" propri\'et\'es.
\vskip 2mm
\noindent
{\bf Th\'eor\`eme\num\label{th2brubru}} {\it 
{\rm 1)} Si $\mathscr C$ est une formation, alors : 
\vskip 2mm
\noindent
$\bullet$ l'extension $k^{\tiny \mathscr C}/k$ est l'unique extension galoisienne $M/k$ v\'erifiant la propri\'et\'e suivante : pour toute extension galoisienne finie $L/k$, on a l'\'equivalence $\hbox{\rm Gal}(L/k)\in {\mathscr C}\Longleftrightarrow L\subseteq M$,
\vskip 2mm
\noindent
$\bullet$ l'extension $k^{\tiny \mathscr C}/k$ est la plus grande extension galoisienne de $k$ ayant pour groupe de Galois un pro-$\mathscr C$-groupe,
\vskip 2mm
\noindent
$\bullet$ si $G\in {\mathscr C}^*$ est le groupe de Galois d'une extension galoisienne finie $E/k$, alors, pour tout corps interm\'ediaire $k\subseteq L\subseteq k^{\tiny \mathscr C}$, le groupe $G$ est aussi le groupe de Galois de l'extension $E \cdot L/L$.
\vskip 2mm
\noindent
{\rm 2)} Si $\mathscr C$ est une formation extensive, alors le corps $k^{\tiny \mathscr C}$ est ${\mathscr C}$-clos. En particulier, si un groupe fini $G$ donn\'e est groupe de Galois sur $k^{\tiny \mathscr C}$, alors $G\in {\mathscr C}^*$. On a donc $(k^{\tiny {\mathscr C}})^{{\mathscr C}^*}=k^{\hbox{\rm \scriptsize s\'ep}}$.
\vskip 2mm
\noindent
{\rm 3)} Si $\mathscr C$ est une vari\'et\'e extensive, alors le corps $k^{\tiny \mathscr C}$ est le plus petit corps contenant $k$ qui soit ${\mathscr C}$-clos. Dans ce cas, si $k_1\subseteq k_2$ d\'esignent deux corps quelconques, alors $k_1^{\tiny \mathscr C}\subseteq k_2^{\tiny \mathscr C}$ (en particulier, si $k_1\subseteq k_2\subseteq k_1^{\tiny \mathscr C}$, alors $k_1^{\tiny \mathscr C}= k_2^{\tiny \mathscr C}$).}
\vskip 2mm
\noindent
{\bf Preuve :} 1) $\bullet$ Consid\'erons une extension galoisienne finie $L/k$ quelconque. Si $\hbox{\rm Gal}(L/k)\in {\mathscr C}$, alors on a $L\subseteq k^{\tiny \mathscr C}$ par d\'efinition de $k^{\mathscr C}$. R\'eciproquement, si $L\subseteq k^{\tiny \mathscr C}$, alors, comme $[L:k]$ est fini, il existe des extensions galoisiennes $E_1/k,\dots ,E_n/k$, \`a groupes de Galois \'el\'ements de ${\mathscr C}$, telles que $L$ soit inclus dans le compositum $E=E_1 \cdot \cdots \cdot E_n$. L'extension finie $E/k$ \'etant galoisienne et son groupe de Galois \'etant un produit fibr\'e d'un nombre fini d'\'el\'ements de ${\mathscr C}$, on a $\hbox{\rm Gal}(E/k)\in {\mathscr C}$ par $(\hbox{\rm C}_3)$. Le groupe $\hbox{\rm Gal}(L/k)$ \'etant un quotient de $\hbox{\rm Gal}(E/k)$, il est \'el\'ement de $\mathscr C$ par $(\hbox{\rm C}_1)$.
\vskip 2mm
\noindent
\tab L'unicit\'e du corps $k^{\mathscr C}$, pour cette propri\'et\'e, est alors \'evidente puisque toute extension galoisienne est la r\'eunion de ses sous-extensions galoisiennes finies.
\vskip 2mm
\noindent
$\bullet$ Les groupes de Galois des extensions galoisiennes finies de $k$ incluses dans $k^{\tiny \mathscr C}$ sont dans ${\mathscr C}$ par ce qui pr\'ec\`ede. Par cons\'equent, $\hbox{\rm Gal}(k^{\tiny \mathscr C}/k)$ est un pro-${\mathscr C}$-groupe. On se donne maintenant une extension galoisienne $M/k$ \`a groupe de Galois pro-${\mathscr C}$. Par d\'efinition, il existe un syst\`eme projectif filtrant \`a droite $(G_i)_i$ d'\'el\'ements de ${\mathscr C}$ tels que $G=\hbox{\rm Gal}(M/k)=\varprojlim_i G_i$. La famille $\{U_i\}_i$ des sous-groupes ouverts normaux associ\'es aux $G_i$ (pour tout $i$, $G/U_i \cong G_i$ pour les projections canoniques) constitue alors un syst\`eme fondamental de voisinages du neutre de $G$. \'Etant donn\'e $\alpha \in M$, notons $\widetilde{k(\alpha)}$ la cl\^oture galoisienne de $k(\alpha)$ sur $k$. Puisque $\widetilde{k(\alpha)}\subseteq M$, le sous-groupe $\hbox{\rm Gal}(M/\widetilde{k(\alpha)})$ de $G$ est un sous-groupe ouvert. Il contient donc un $U_i$ pour un certain $i$, et il s'ensuit que $\hbox{\rm Gal}(\widetilde{k(\alpha)}/k)$ est un quotient du groupe $G_i$. Comme $G_i \in {\mathscr C}$ et ${\mathscr C}$ v\'erifie $(\hbox{\rm C}_1)$, l'on voit que $\hbox{\rm Gal}(\widetilde{k(\alpha)}/k)$ est \'el\'ement de ${\mathscr C}$. Ainsi $\alpha \in \widetilde{k(\alpha)}\subseteq k^{\tiny \mathscr C}$.
\vskip 2mm
\noindent
$\bullet$ Le groupe de Galois de l'extension galoisienne finie $E\cap k^{\tiny \mathscr C}/k$ est, d'une part, un quotient de $G$, donc dans ${\mathscr C}^*$ (en vertu du 1) de la proposition \ref{transmission 2}), et, d'autre part, un \'el\'ement de ${\mathscr C}$ d'apr\`es le premier point ci-dessus. Il est donc trivial. Ainsi on a $E\cap k^{\tiny \mathscr C}=k$, ce qui \'equivaut \`a dire que les corps $E$ et $k^{\tiny \mathscr C}$ sont lin\'eairement disjoints sur $k$. On en d\'eduit que, pour tout corps interm\'ediaire $k\subseteq L\subseteq k^{\tiny \mathscr C}$, les corps $E$ et $L$ sont \'egalement lin\'eairement disjoints sur $k$. Par cons\'equent, on a $\hbox{\rm Gal}(E \cdot L/L)=G$.
\vskip 2mm
\noindent
2) Consid\'erons une extension galoisienne finie $L/k^{\tiny \mathscr C}$ telle que $G=\hbox{\rm Gal}(L/k^{\tiny \mathscr C})$ soit \'el\'ement de ${\mathscr C}$ et un \'el\'ement $\alpha$ de $L$ tel que $L=k^{\tiny \mathscr C}(\alpha)$.
\vskip 2mm
\noindent
\tab D\'emontrons tout d'abord que $L=k^{\mathscr C}$ sous l'hypoth\`ese suppl\'ementaire que l'extension $L/k$ soit galoisienne \`a groupe de Galois \'el\'ement de ${\mathscr C}$. Pour cela, notons $P(X) \in k^{\mathscr C}[X]$ le polyn\^ome minimal de $\alpha$ sur $k^{\tiny \mathscr C}$. Comme l'extension $L/k$ est galoisienne, pour tout $\sigma \in \hbox{\rm Gal}(k^{\hbox{\scriptsize s\'ep}}/k)$, il existe $P_\sigma(X) \in k^{\tiny \mathscr C}[X]$ tel que $\sigma(\alpha)=P_\sigma(\alpha)$. Puisque $\alpha$ est alg\'ebrique sur $k$, l'on peut choisir $\sigma_1, \dots, \sigma_n \in \hbox{\rm Gal}(k^{\hbox{\scriptsize s\'ep}}/k)$ tels que $P_\sigma(X) \in \{P_{\sigma_1}(X), \dots, P_{\sigma_n}(X)\}$ pour tout $\sigma \in \hbox{\rm Gal}(k^{\hbox{\scriptsize s\'ep}}/k)$. Consid\'erons alors le corps $\Omega$, cl\^oture galoisienne sur $k$ du corps obtenu en adjoignant \`a $k$ les coefficients des polyn\^omes $P_{\sigma_1}(X), \dots, P_{\sigma_n}(X)$. Par construction, l'extension $\Omega(\alpha)/k$ est galoisienne, et l'on a $\hbox{\rm Gal}(\Omega(\alpha)/\Omega)=G$. Puisque $\Omega/k$ est une extension galoisienne finie incluse dans $k^{\tiny \mathscr C}$, on a $\hbox{\rm Gal}(\Omega/k)\in {\mathscr C}$ d'apr\`es le premier point du 1) ci-dessus. La classe ${\mathscr C}$ v\'erifiant $(\hbox{\rm C}_2)$, l'on voit alors que le groupe $\hbox{\rm Gal}(\Omega(\alpha)/k)$, en tant qu'extension de $\hbox{\rm Gal}(\Omega(\alpha)/\Omega)=G\in {\mathscr C}$ par $\hbox{\rm Gal}(\Omega/k)\in {\mathscr C}$, est aussi un \'el\'ement de ${\mathscr C}$. En particulier, on a $\alpha \in \Omega(\alpha)\subseteq k^{\tiny \mathscr C}$ et donc $L=k^{\tiny \mathscr C}$.
\vskip 2mm
\noindent
\tab Nous expliquons maintenant comment traiter le cas g\'en\'eral. Pour tout $\sigma \in \hbox{\rm Gal}(k^{\hbox{\scriptsize s\'ep}}/k)$, l'extension finie $k^{\tiny \mathscr C}(\sigma(\alpha))/k^{\tiny \mathscr C}$ est galoisienne de groupe de Galois $G$. En effet, dans $\hbox{\rm Gal}(k^{\hbox{\scriptsize s\'ep}}/k)$, le sous-groupe $\hbox{\rm Gal}(k^{\hbox{\scriptsize s\'ep}}/k^{\tiny \mathscr C}(\sigma(\alpha)))$ est le conjugu\'e par $\sigma$ du sous-groupe $\hbox{\rm Gal}(k^{\hbox{\scriptsize s\'ep}}/k^{\tiny \mathscr C}(\alpha))$. Puisque la conjugaison par $\sigma$ induit un automorphisme de $ \hbox{\rm Gal}(k^{\hbox{\scriptsize s\'ep}}/k^{\tiny \mathscr C})$, on a alors $\hbox{\rm Gal}(k^{\hbox{\scriptsize s\'ep}}/k^{\tiny \mathscr C}(\sigma(\alpha)))\trianglelefteq \hbox{\rm Gal}(k^{\hbox{\scriptsize s\'ep}}/k^{\tiny \mathscr C})$ et 
$$\begin{array}{lll}
\hbox{\rm Gal}(k^{\tiny \mathscr C}(\sigma(\alpha))/k^{\tiny \mathscr C})&\cong &\hbox{\rm Gal}(k^{\hbox{\scriptsize s\'ep}}/k^{\tiny \mathscr C})/\hbox{\rm Gal}(k^{\hbox{\scriptsize s\'ep}}/k^{\tiny \mathscr C}(\sigma(\alpha)))\\
&=&\hbox{\rm Gal}(k^{\hbox{\scriptsize s\'ep}}/k^{\tiny \mathscr C})/\sigma\hbox{\rm Gal}(k^{\hbox{\scriptsize s\'ep}}/k^{\tiny \mathscr C}(\alpha))\sigma^{-1}\\
&\cong &\hbox{\rm Gal}(k^{\hbox{\scriptsize s\'ep}}/k^{\tiny \mathscr C})/\hbox{\rm Gal}(k^{\hbox{\scriptsize s\'ep}}/k^{\tiny \mathscr C}(\alpha))\\
&\cong& {\rm Gal}(k^{\mathscr C}(\alpha)/k^{\mathscr C})\\
&= & G.\\
\end{array}$$
Le compositum $M=\bullet_{\sigma} k^{\tiny \mathscr C}(\sigma(\alpha))$ d\'efinit donc une extension galoisienne de $k$ (c'est en fait la cl\^oture galoisienne de $k^{\tiny \mathscr C}(\alpha)$ sur $k$) et le groupe de Galois $\hbox{\rm Gal}(M/k^{\tiny \mathscr C})$ est alors le produit fibr\'e d'un nombre fini de copies du groupe $G$. Puisque ${\mathscr C}$ v\'erifie $(\hbox{\rm C}_3)$, ce groupe est dans ${\mathscr C}$. Par ce qui pr\'ec\`ede, on a alors $M=k^{\mathscr C}$ et donc $L=k^{\mathscr C}$.
\vskip 2mm
\noindent
\tab On se donne enfin une extension galoisienne finie $L/k^{\tiny \mathscr C}$ de groupe de Galois $G$ et un sous-groupe normal $H$ de $G$. Alors la sous-extension $L^H/k^{\tiny \mathscr C}$  de $L/k^{\mathscr C}$ est galoisienne finie de groupe de Galois $G/H$. Si $G/H\in {\mathscr C}$, alors $L^H=k^{\tiny \mathscr C}$ d'apr\`es ce qui pr\'ec\`ede, et donc $H=G$. Ainsi $G\in {\mathscr C}^*$.
\vskip 2mm
\noindent
3) On sait d\'ej\`a, d'apr\`es 2), que $k^{\tiny \mathscr C}$ est ${\mathscr C}$-clos. On se donne maintenant un corps $L$ contenant $k$ qui soit ${\mathscr C}$-clos et une extension galoisienne finie $M/k$ \`a groupe de Galois dans ${\mathscr C}$. Comme le groupe de Galois de l'extension galoisienne finie $M \cdot L/L$ s'identifie \`a un sous-groupe de $\hbox{\rm Gal}(M/k)\in {\mathscr C}$, le groupe ${\rm{Gal}}(M \cdot L/L)$ est lui aussi \'el\'ement de ${\mathscr C}$ (par $(\hbox{\rm C}_0)$). Le corps $L$ \'etant ${\mathscr C}$-clos, on a alors $M \cdot L=L$, c'est-\`a-dire $M\subseteq L$. On a donc $k^{\tiny \mathscr C}\subseteq L$.
\vskip 2mm
\noindent
\tab Si $k_1\subseteq k_2$ d\'esignent deux corps quelconques, alors $k_1\subseteq k_2\subseteq k_2^{\tiny \mathscr C}$. Puisque $k_2^{\tiny \mathscr C}$ est ${\mathscr C}$-clos par ce qui pr\'ec\`ede, on a $k_1^{\tiny \mathscr C}\subseteq k_2^{\tiny \mathscr C}$ par minimalit\'e de $k_1^{\tiny \mathscr C}$.
\fin
\vskip 2mm
\noindent
{\bf Remarque :} Si ${\mathscr C}$ est une formation extensive, on voit que, d'apr\`es le 2) du th\'eor\`eme \ref{th2brubru}, aucun \'el\'ement non trivial de ${\mathscr C}$ ne se r\'ealise comme groupe de Galois sur $k^{\tiny \mathscr C}$. Si ${\mathscr C}$ contient de plus un groupe de Galois r\'egulier sur $k$ (par exemple, un groupe ab\'elien fini ou un groupe sym\'etrique), l'on voit que le corps $k^{\mathscr C}$ n'est pas hilbertien. Il est \'evident que
$${\mathscr C}={\mathscr C}(1)\ \hbox{\rm et}\ k\ \hbox{\rm hilbertien} \Longrightarrow k^{\mathscr C}\ \hbox{\rm  hilbertien}.$$
\vskip 2mm
\noindent
La r\'eciproque de cette implication est vraie si le PIGR$_{/k}$ admet une r\'eponse positive, en vertu de ce qui pr\'ec\`ede. Elle est aussi vraie si ${\mathscr C}$ v\'erifie ${\rm{(C_0)}}$. En effet, si $k^{\mathscr C}$ est hilbertien, le raisonnement ci-dessus montre alors que ${\mathscr C}$ ne contient aucun groupe ab\'elien fini non trivial. Puisque ${\mathscr C}$ v\'erifie ${\rm{(C_0)}}$, le th\'eor\`eme de Cauchy assure alors que ${\mathscr C}={\mathscr C}$(1). Ainsi $k^{\mathscr C}=k$ et $k$ est hilbertien.
\vskip 2mm
\noindent
\tab Nous donnons ci-dessous une variante du troisi\`eme point du 1) du th\'eor\`eme \ref{th2brubru} valable sous la seule hypoth\`ese que $\mathscr{C}$ v\'erifie $({\rm{C}}_1)$.
\vskip 2mm
\noindent
{\bf{Proposition\num\label{propbrubru}}} {\it{Si $\mathscr C$ est une pr\'e-formation et si un groupe fini $G\in {\mathscr C}^*$ n'ayant aucun quotient d'ordre premier est le groupe de Galois d'une extension galoisienne finie $E/k$, alors, pour tout corps interm\'ediaire $k\subseteq L\subseteq k^{\mathscr C}$, le groupe $G$ est aussi le groupe de Galois de l'extension $E \cdot L/L$.}}
\vskip 2mm
\noindent
{\bf Preuve :} D\'emontrons tout d'abord un r\'esultat pr\'eliminaire sur les quotients simples d'un produit fibr\'e de deux groupes quelconques :
\vskip 2mm
\noindent
{\it{\'Etant donn\'es un groupe simple non ab\'elien $S$, un groupe $G$ et deux sous-groupes normaux $H_1$ et $H_2$ de $G$, si $S$ n'est quotient ni de $G/H_1$, ni de $G/H_2$, alors $S$ n'est quotient de $G/(H_1 \cap H_2)$.}}
\vskip2mm
\noindent
En effet, quitte \`a quotienter $G$, $H_1$ et $H_2$ par $H_1 \cap H_2$, on peut supposer que $H_1$ et $H_2$ sont d'intersection triviale. S'il existe un sous-groupe normal $N$ de $G$ tel que $G/N=S$, alors on a n\'ecessairement $H_1N=H_2N=G$ puisque le groupe simple $S$ n'est quotient ni de $G/H_1$, ni de $G/H_2$. L'observation suivante, qui nous a \'et\'e communiqu\'ee par A. Fehm et due \`a M. Shusterman, permet alors de conclure :
\vskip 2mm
\noindent
{\it Si $H$ d\'esigne un groupe et $U,V$ deux sous-groupes normaux de $H$ d'intersection triviale, alors, pour tout sous-groupe normal $K$ de $H$ tel que $KU=KV=H$, le groupe $H/K$ est ab\'elien\footnote{En effet, si $x$ et $y$ sont deux \'el\'ements quelconques de $H$, on a $x=k_1u$ et $u=k_2v$ avec $u\in U$, $v\in V$ et $(k_1,k_2)^2\in K$. On a donc
$[x,y]=k_1uk_2vu^{-1}k_1^{-1}v^{-1}k_2^{-1}=k_1(uk_2u^{-1})uvu^{-1}v^{-1}(vk_1^{-1}v^{-1})k_2^{-1}.$
Il ne reste plus qu'\`a remarquer que $[U,V]=\{1\}$ (puisque $U$ et $V$ sont normaux et d'intersection triviale) pour conclure que $[x,y]$ est \'el\'ement de $K$.}.}
\vskip 2mm
\noindent
\tab Venons-en maintenant \`a la d\'emonstration de la proposition. Clairement, il suffit de d\'emontrer que $E \cap k^{\tiny \mathscr C}=k$. Si le groupe $\hbox{\rm Gal}((E \cap k^{\tiny \mathscr C})/ k)$ n'est pas trivial, alors ce groupe poss\`ede un quotient simple $S$ qui, en tant que quotient de $G$, est non ab\'elien et \'el\'ement de ${\mathscr C}^*$ (d'apr\`es le 1) de la proposition \ref{transmission 2}), et n'est donc pas \'el\'ement de ${\mathscr C}$. Si $L_0$ d\'esigne le sous-corps de $E \cap k^{\tiny \mathscr C}$ tel que ${\rm Gal}(L_0/k)=S$, alors on a  $L_0 \subseteq k^{\tiny \mathscr C}$. Puisque $[L_0:k]$ est fini, il existe des extensions galoisiennes finies $M_1/k,\dots ,M_n/k$, \`a groupes de Galois \'el\'ements de ${\mathscr C}$, telles que $L_0$ soit inclus dans le compositum $M=M_1 \cdot \cdots \cdot M_n$. La classe $\mathscr{C}$ v\'erifiant $({\rm{C}}_1)$, le groupe $S$ n'est quotient de ${\rm{Gal}}(M_i/k)$ pour aucun $i \in \{1,\dots,n\}$. En vertu du r\'esultat pr\'eliminaire ci-dessus et d'une r\'ecurrence imm\'ediate, le groupe ${\rm{Gal}}(M/k)$, qui s'identifie \`a un produit fibr\'e successif des groupes ${\rm{Gal}}(M_1/k), \dots, {\rm{Gal}}(M_n/k)$, n'admet alors pas $S$ comme quotient, ce qui est manifestement impossible puisque $L_0$ est un sous-corps de $M$. Par cons\'equent, $\hbox{\rm Gal}((E \cap k^{\tiny \mathscr C})/ k)$ est trivial, c'est-\`a-dire $E \cap k^{\tiny \mathscr C}=k$.
\fin
\vskip 2mm
\noindent
{\bf Remarque :} Si $\mathscr{C}$ contient tous les groupes finis d'ordre premier, la condition qu'aucun quotient simple de $G$ ne soit ab\'elien dans la proposition \ref{propbrubru} est automatique si $G \in \mathscr{C}^*$  (par le 1) de la proposition \ref{transmission 2}).
\vskip 2mm
\noindent
\tab On d\'efinit maintenant la suite $(k^{\tiny {\mathscr C}[n]})_{n \geq 0}$ des {\it ${\mathscr C}$-cl\^otures successives} du corps $k$, en posant $k^{\tiny {\mathscr C}[0]}=k$ et $k^{\tiny {\mathscr C}[n+1]}=(k^{\tiny {\mathscr C}[n]})^{\tiny \mathscr C}$ pour tout $n\geq 0$. 
\vskip 2mm
\noindent
\tab La r\'eunion des ${\mathscr C}$-cl\^otures successives de $k$ est alors reli\'ee \`a la  $\widehat{\mathscr C}$-cl\^oture de $k$, de la mani\`ere suivante :
\vskip 2mm
\noindent
{\bf Th\'eor\`eme \num\label{th3brubru}} {\it {\rm 1)} Pour tout $n \geq 0$, l'extension $k^{\tiny {\mathscr C}[n]}/k$ est galoisienne.
\vskip 2mm
\noindent
{\rm 2)} Si ${\mathscr C}$ est une pr\'e-vari\'et\'e, alors $\displaystyle k^{\tiny \widehat{\mathscr C}} \subseteq \bigcup_{n\geq 0}k^{\tiny {\mathscr C}[n]}$.
\vskip 2mm
\noindent
{\rm 3)} Si ${\mathscr C}$ est une pr\'e-vari\'et\'e qui v\'erifie en plus la condition $(\hbox{\rm C}_3)$, alors $\displaystyle k^{\tiny \widehat{\mathscr C}} = \bigcup_{n\geq 0}k^{\tiny {\mathscr C}[n]}$.}
\vskip 2mm
\noindent
{\bf Preuve :} 1) Pour $n \in \{0,1\}$, l'extension $k^{\tiny {\mathscr C}[n]}/k$ est galoisienne. On se donne maintenant $n \geq 1$ et l'on suppose que $k^{\tiny {\mathscr C}[n]}/k$ est galoisienne. \'Etant donn\'e $\alpha \in k^{\tiny {\mathscr C}[n+1]}$, il existe des extensions galoisiennes finies $M_1/k^{\tiny {\mathscr C}[n]}, \dots, M_m/k^{\tiny {\mathscr C}[n]}$, \`a groupes de Galois dans ${\mathscr C}$, telles que $\alpha$ appartienne au compositum $M=M_1 \cdot \cdots \cdot M_m$. Fixons alors $i \in \{1,\dots,m\}$, un \'el\'ement primitif $\beta_i$ de $M_i/k^{\tiny {\mathscr C}[n]}$ et $\sigma \in \hbox{\rm Gal}(k^{\hbox{\scriptsize s\'ep}}/k)$. En utilisant uniquement l'hypoth\`ese $k^{\tiny {\mathscr C}[n]}/k$ galoisienne, l'on montre comme dans la preuve du 2) du th\'eor\`eme \ref{th2brubru} que l'extension finie $k^{\tiny {\mathscr C}[n]}(\sigma(\beta_i))/k^{\tiny {\mathscr C}[n]}$ est galoisienne de groupe de Galois $G_i ={\rm Gal}(M_i/k^{\tiny {\mathscr C}[n]})$. Comme $G_i \in {\mathscr C}$, l'on voit que $\sigma(\beta_i) \in k^{\tiny {\mathscr C}[n]}(\sigma(\beta_i)) \subseteq k^{\tiny {\mathscr C}[n+1]}$, ce qui montre que $\sigma(M_i) \subseteq k^{\tiny {\mathscr C}[n+1]}$. En particulier, on a $\sigma(\alpha) \in \sigma(M) \subseteq k^{\tiny {\mathscr C}[n+1]}$.
\vskip 2mm
\noindent
2) \'Etant donn\'e $\alpha \in k^{\tiny \widehat{\mathscr C}}$, consid\'erons la cl\^oture galoisienne $M$ de $k(\alpha)$ sur $k$; celle-ci est contenue dans $k^{\tiny \widehat{\mathscr C}}$ puisque $k^{\tiny \widehat{\mathscr C}}/k$ est galoisienne. Or $\widehat{{\mathscr C}}$ est une formation puisque ${\mathscr C}$ est une pr\'e-vari\'et\'e (cf. corollaire \ref{ffr}). Par le premier point du 1) du th\'eor\`eme \ref{th2brubru}, on a donc $G={\rm Gal}(M/k) \in \widehat{\mathscr C}$. Par d\'efinition de $\widehat{{\mathscr C}}$, il existe alors une tour d'extensions finies $k=M_0\subseteq M_1\subseteq \cdots \subseteq M_{n-1} \subseteq M_n=M$ telle que, pour tout $i \in \{0, \dots, n-1\}$, l'extension $M_{i+1}/M_i$ soit galoisienne de groupe de Galois $G_i \in {\mathscr C}$. Pour tout $i \in \{0, \dots ,n-1\}$, le groupe de Galois de l'extension galoisienne finie $M_{i+1} \cdot k^{\tiny {\mathscr C}[i]} / M_i \cdot k^{\tiny {\mathscr C}[i]}$ s'identifie \`a un sous-groupe de $G_i$. Or $G_i$ est dans ${\mathscr C}$ et ${\mathscr C}$ v\'erifie $(\hbox{\rm C}_0)$. On a donc ${\rm Gal}(M_{i+1} \cdot k^{\tiny {\mathscr C}[i]} / M_i \cdot k^{\tiny {\mathscr C}[i]}) \in {\mathscr C}$, et une r\'ecurrence imm\'ediate montre alors que $M=M_n \subseteq k^{\tiny {\mathscr C}[n]}$. On a donc $\alpha \in M \subseteq k^{\tiny {\mathscr C}[n]}$.
\vskip 1mm
\noindent
3) R\'eciproquement, on a $k^{\tiny {\mathscr C}[0]}=k\subseteq k^{\tiny \widehat{\mathscr C}}$. On se donne maintenant un entier $n \geq 0$ et l'on suppose que $k^{\tiny {\mathscr C}[n]}\subseteq k^{\tiny \widehat{\mathscr C}}$. \'Etant donn\'e $\alpha \in k^{\tiny {\mathscr C}[n+1]}$, consid\'erons la cl\^oture galoisienne $M$ de $k^{\tiny {\mathscr C}[n]}(\alpha)$ sur $k^{\tiny {\mathscr C}[n]}$. Comme $M \subseteq k^{\tiny {\mathscr C}[n+1]}$ et ${\mathscr C}$ est une formation, on peut \`a nouveau utiliser le premier point du 1) du th\'eor\`eme \ref{th2brubru} : on a ${\rm{Gal}}(M/k^{\tiny {\mathscr C}[n]}) \in {\mathscr C}$. Or, par hypoth\`ese, on a $k^{\tiny {\mathscr C}[n]} \subseteq k^{\tiny \widehat{\mathscr C}}$. Par cons\'equent, ${\rm{Gal}}(M \cdot k^{\tiny \widehat{\mathscr C}} / k^{\tiny \widehat{\mathscr C}})$ s'identifie \`a un sous-groupe de ${\rm Gal}(M/k^{\tiny {\mathscr C}[n]}) \in {\mathscr C}$. Puisque ${\mathscr C}$ v\'erifie $(\hbox{\rm C}_0)$ et ${\mathscr C} \subseteq \widehat{{\mathscr C}}$, l'on voit que ${\rm{Gal}}(M \cdot k^{\tiny \widehat{\mathscr C}} / k^{\tiny \widehat{\mathscr C}}) \in \widehat{{\mathscr C}}$. La classe ${\mathscr C}$ \'etant une pr\'e-vari\'et\'e, la classe $\widehat{{\mathscr C}}$ est une formation extensive (cf. corollaire \ref{ffr}). En vertu du 2) du th\'eor\`eme \ref{th2brubru}, on a alors $M \cdot k^{\tiny \widehat{\mathscr C}}=k^{\tiny \widehat{\mathscr C}}$, c'est-\`a-dire $M \subseteq k^{\tiny \widehat{\mathscr C}}$. En particulier, on a $\alpha \in M \subseteq k^{\tiny \widehat{\mathscr C}}$.
\fin
\vskip 5mm 
\noindent
{\bf 4.2.--- Groupes de Galois sur une ${\mathscr C}$-cl\^oture.}
\vskip 2mm 
\noindent
\tab Le 2) du th\'eor\`eme \ref{th2brubru} assure que, si ${\mathscr C}$ est une formation extensive, alors les groupes finis qui apparaissent comme groupes de Galois sur le corps $k^{\tiny \mathscr C}$ sont n\'ecessairement \'el\'ements de ${\mathscr C}^*$. La r\'eciproque de cette propri\'et\'e est fausse en g\'en\'eral puisque, par exemple, la cl\^oture r\'esoluble de ${\mathbb F}_p$ est \'egale \`a $\overline{\mathbb F}_p$ pour tout nombre premier $p$. Par application du troisi\`eme point du 1) du th\'eor\`eme \ref{th2brubru}, on peut toutefois remarquer que, si le $\hbox{\rm PIG}_{/k}$ admet une r\'eponse positive, alors tous les groupes finis appartenant \`a ${\mathscr C}^*$ sont bien groupes de Galois sur $ k^{\tiny \mathscr C}$. Il est donc assez raisonnable de conjecturer que, si $k$ est un corps hilbertien, alors les groupes finis apparaissant comme groupes de Galois sur $ k^{\tiny \mathscr C}$ sont exactement les \'el\'ements de ${\mathscr C}^*$.
\vskip 2mm
\noindent
\tab Nous allons maintenant d\'ecrire une autre situation dans laquelle les groupes finis apparaissant comme groupes de Galois sur $k^{\tiny \mathscr C}$ sont exactement les \'el\'ements de ${\mathscr C}^*$, et faire ensuite le lien avec une conjecture tr\`es c\'el\`ebre de th\'eorie inverse de Galois\footnote{La proposition \ref{sha} est librement inspir\'ee d'une id\'ee de D. Haran et M. Jarden communiqu\'ee au premier auteur de cet article dans une correspondance au sujet de la cl\^oture r\'esoluble de $\mathbb Q$.}.
\vskip 2mm
\noindent
{\bf Proposition\num\label{sha}} {\it On suppose que ${\mathscr C}$ est une vari\'et\'e extensive. S'il existe un corps interm\'ediaire $k \subseteq L \subseteq k^{\mathscr C}$ \`a groupe de Galois absolu prolibre de rang $\alpha$ infini, alors les groupes profinis qui apparaissent comme groupes de Galois sur $k^{\tiny \mathscr C}$ sont exactement les pro-${\mathscr C}^*$-groupes de rang $\leq \alpha$. 
\vskip 2mm
\noindent
En particulier, dans cette situation, les groupes finis qui apparaissent comme groupes de Galois sur $ k^{\tiny \mathscr C}$ sont exactement les \'el\'ements de ${\mathscr C}^*$.}
\vskip 2mm
\noindent
{\bf Preuve :} \'Etant donn\'e un corps interm\'ediaire $L$ comme ci-dessus, commen\c cons par remarquer que, d'apr\`es le 3) du th\'eor\`eme \ref{th2brubru}, on a $L^{\mathscr C}=k^{\tiny \mathscr C}$. Ainsi, si l'on identifie le groupe de Galois absolu de $L$, $\hbox{\rm Gal}(k^{\hbox{\scriptsize s\'ep}}/L)$, au groupe prolibre $\widehat{F}_\alpha$, on voit que le groupe $N=\hbox{\rm Gal}(k^{\hbox{\scriptsize s\'ep}}/ k^{\tiny \mathscr C})$ (qui, d'apr\`es le premier point du 1) du th\'eor\`eme \ref{th2brubru}, est \'egal \`a l'intersection des sous-groupes ouverts normaux $\Gamma$ de $\widehat{F}_\alpha$ tels que $\widehat{F}_\alpha/\Gamma$ soit \'el\'ement de ${\mathscr C}$) v\'erifie les deux propri\'et\'es suivantes :
\vskip 1mm
\noindent
$\bullet$ {\hbox{$N$ est de rang $\leq \alpha$ (en tant que sous-groupe ferm\'e d'un groupe profini de rang infini $\alpha$),}}
\vskip 1mm
\noindent
$\bullet$ $N$ est un pro-${\mathscr C}^*$-groupe (car les quotients finis de $N$ correspondent aux groupes de Galois des extensions galoisiennes finies de $k^{\tiny \mathscr C}$ qui, d'apr\`es le 2) du th\'eor\`eme \ref{th2brubru}, sont tous \'el\'ements de ${\mathscr C}^*$).
\vskip 2mm
\noindent
Les groupes profinis qui sont groupes de Galois sur $k^{\tiny \mathscr C}$ sont exactement les quotients de $N$. Ils ont donc tous un rang $\leq \alpha$ et sont des pro-${\mathscr C}^*$-groupes car ${\mathscr C}^*$ v\'erifie $(\hbox{\rm C}_1)$ (par le 1) de la proposition \ref{transmission 2}).
\vskip 2mm
\noindent
\tab R\'eciproquement, l'on se donne un pro-${\mathscr C}^*$-groupe $G$ de rang $\leq \alpha$. Par propri\'et\'e universelle de $\widehat{F}_\alpha$, il existe un \'epimorphisme $\varphi:\widehat{F}_\alpha\longrightarrow G$ qui induit alors, par passage au quotient, un \'epimorphisme $\widetilde{\varphi}:\widehat{F}_\alpha/N\longrightarrow G/\varphi(N)$. Par construction, le groupe quotient $\widehat{F}_\alpha/N$ est le groupe de Galois de l'extension galoisienne $L^{\tiny \mathscr C}/L$;  c'est donc un pro-${\mathscr C}$-groupe en vertu du deuxi\`eme point du 1) du th\'eor\`eme \ref{th2brubru}. Puisque ${\mathscr C}$ v\'erifie $(\hbox{\rm C}_1)$, l'on voit que, via l'\'epimorphisme $\widetilde{\varphi}$, le groupe quotient $G/\varphi(N)$ est lui aussi un pro-${\mathscr C}$-groupe. Tout quotient fini non trivial de $G/\varphi(N)$ est alors un \'el\'ement de ${\mathscr C}$, mais aussi un \'el\'ement de ${\mathscr C}^*$ en tant que quotient du pro-${\mathscr C}^*$-groupe $G$. On en d\'eduit que $G/\varphi(N)$ n'a pas de quotient non trivial et donc que $G=\varphi(N)$. Ainsi $G$ est un quotient de $N$ et est donc groupe de Galois d'une extension galoisienne de $k^{\tiny \mathscr C}$. 
\fin
\vskip 2mm
\noindent
\tab Pour illustrer l'int\'er\^et de la proposition \ref{sha}, nous consid\'erons la classe ${\mathscr C}(\hbox{\rm r\'es})$ des groupes finis r\'esolubles. Comme c'est une vari\'et\'e extensive, les th\'eor\`emes \ref{th2brubru} et \ref{th3brubru}, ainsi que la proposition \ref{sha}, peuvent alors \^etre enti\`erement appliqu\'es. En particulier, la cl\^oture r\'esoluble de $k$, not\'ee $k^{\hbox{\rm \scriptsize r\'es}}$, est l'unique corps $M$ contenant $k$ v\'erifiant les conditions (\'equivalentes) suivantes :
\vskip 1mm
\noindent
i) l'extension $M/k$ est galoisienne et, pour toute extension galoisienne finie $L/k$, le groupe $\hbox{\rm Gal}(L/k)$ est r\'esoluble si et seulement si $L\subseteq M$,
\vskip 1mm
\noindent
ii) l'extension $M/k$ est la plus grande extension galoisienne de $k$ ayant pour groupe de Galois un groupe pro-r\'esoluble,
\vskip 1mm
\noindent
iii) le corps $M$ est le plus petit corps contenant $k$ qui soit ${\mathscr C}(\hbox{\rm r\'es})$-clos,
\vskip 1mm
\noindent
iv) le corps $M$ est \'egal \`a la r\'eunion des ab\'elianis\'es successifs de $k$.
\vskip 2mm
\noindent
\tab Puisque l'on a les inclusions ${\mathscr C}(\hbox{\rm cycl})\subseteq {\mathscr C}(\hbox{\rm ab})\subseteq {\mathscr C}(\hbox{\rm nil})\subseteq {\mathscr C}(\hbox{\rm r\'es})$, on en d\'eduit que ${\mathscr C}(\hbox{\rm cycl})^*\supseteq {\mathscr C}(\hbox{\rm ab})^*\supseteq {\mathscr C}(\hbox{\rm nil})^*\supseteq {\mathscr C}(\hbox{\rm r\'es})^*$. Le 3) de la proposition \ref{transmission} montre que ces inclusions sont en fait des \'egalit\'es.
\vskip 2mm
\noindent
\tab Par d\'efinition, la classe des groupes finis fortement non r\'esolubles, not\'ee ${\mathscr C}(\hbox{\rm FnR})$, est la classe duale ${\mathscr C}(\hbox{\rm cycl})^*={\mathscr C}(\hbox{\rm ab})^*={\mathscr C}(\hbox{\rm nil})^*={\mathscr C}(\hbox{\rm r\'es})^*$. On a alors les propri\'et\'es suivantes :
\vskip 2mm
\noindent
a) Un groupe fini est dans ${\mathscr C}(\hbox{\rm FnR})$ si et seulement si aucun de ses quotients n'est d'ordre premier.
\vskip 2mm
\noindent
b) On a $\displaystyle {\mathscr C}(\hbox{\rm FnR})=\bigcap_{p\in {\mathscr P}}{\mathscr C}(\hbox{\rm p})^*$ o\`u ${\mathscr P}$ d\'esigne l'ensemble des nombres premiers.
\vskip 2mm
\noindent
c) La classe ${\mathscr C}(\hbox{\rm FnR})$ est la plus grande classe de groupes finis v\'erifiant $(\hbox{\rm C}_1)$ et ne contenant aucun groupe fini d'ordre premier.
\vskip 2mm
\noindent
d) La classe ${\mathscr C}(\hbox{\rm FnR})$ v\'erifie $(\hbox{\rm C}_2)$. En particulier, toute extension et tout produit fibr\'e de groupes finis simples non ab\'eliens sont des groupes fortement non r\'esolubles. 
\vskip 2mm
\noindent
e) La classe ${\mathscr C}(\hbox{\rm FnR})$ ne satisfait visiblement pas la condition $(\hbox{\rm C}_0)$. En fait, elle n'est m\^eme pas stable par sous-groupes normaux. En effet, consid\'erons par exemple le groupe altern\'e $A_5$, qui est simple non ab\'elien et donc fortement non r\'esoluble. Puisqu'il se r\'ealise comme groupe de Galois sur ${\mathbb Q}$, il se r\'ealise aussi sur ${\mathbb Q}^{\hbox{\rm \scriptsize r\'es}}$ (en vertu du troisi\`eme point du 1) du th\'eor\`eme \ref{th2brubru}) et l'on peut donc consid\'erer une extension galoisienne finie $L/{\mathbb Q}^{\hbox{\rm \scriptsize r\'es}}$ de groupe de Galois $A_5$. L'extension ${\mathbb Q}^{\hbox{\rm \scriptsize r\'es}}/\mathbb{Q}$ \'etant galoisienne, un c\'el\`ebre th\'eor\`eme de Weissauer\footnote{qui assure que toute extension finie stricte et s\'eparable d'une extension galoisienne d'un corps hilbertien est hilbertienne (voir par exemple \cite[Theorem 13.9.1]{FJ08}).} assure alors que $L$ est hilbertien, et donc que ${\mathbb Z}/2\mathbb{Z}$ se r\'ealise comme groupe de Galois sur $L$. \'Etant donn\'es une extension quadratique $M/L$ et $\sigma \in {\rm Gal}(\overline{\mathbb{Q}}/{\mathbb Q}^{\hbox{\rm \scriptsize r\'es}})$, l'extension $\sigma(M)/L$ est clairement quadratique. Par cons\'equent, si l'on note $\widetilde{M}$ le compositum des corps $\sigma(M)$ lorsque $\sigma$ parcourt le groupe ${\rm Gal}(\overline{\mathbb Q}/{\mathbb Q}^{\hbox{\rm \scriptsize r\'es}})$ ($\widetilde{M}$ est en fait la cl\^oture galoisienne de $M$ sur ${\mathbb Q}^{\hbox{\rm \scriptsize r\'es}}$), l'on voit qu'il existe un entier $N \geq 1$ tel que $\hbox{\rm Gal}(\widetilde{M}/L)=({\mathbb Z}/2\mathbb{Z})^N$. Il s'ensuit que $\hbox{\rm Gal}(\widetilde{M}/{\mathbb Q}^{\hbox{\rm \scriptsize r\'es}})$, qui est dans ${\mathscr C}(\hbox{\rm FnR})$ (puisque se r\'ealisant comme groupe de Galois sur ${\mathbb Q}^{\hbox{\rm \scriptsize r\'es}}$), poss\`ede un sous-groupe normal, $({\mathbb Z}/2\mathbb{Z})^N$, qui, \'etant ab\'elien, ne peut \^etre dans ${\mathscr C}(\hbox{\rm FnR})$.
\vskip2mm
\noindent
f) La classe duale de ${\mathscr C}(\hbox{\rm FnR})$ (c'est-\`a-dire la classe biduale de ${\mathscr C}(\hbox{\rm r\'es})$) est la classe constitu\'ee des groupes finis qui n'admettent pour quotients simples que des groupes d'ordre premier.
\vskip 2mm
\noindent
\tab Regardons maintenant deux applications de la proposition \ref{sha} \`a la classe ${\mathscr C}(\hbox{\rm r\'es})$ :
\vskip 2mm 
\noindent
a) $k={\mathbb Q}$. Comme nous l'avons d\'ej\`a remarqu\'e, une r\'eponse positive au $\hbox{\rm PIG}_{/\mathbb{Q}}$ entra\^inerait que les groupes finis apparaissant comme groupes de Galois sur ${\mathbb Q}^{\hbox{\rm \scriptsize r\'es}}$ seraient exactement les groupes finis fortement non r\'esolubles. Il est int\'eressant de remarquer que, sous la c\'el\`ebre conjecture de Shafarevich\footnote{Cette conjecture affirme que le groupe de Galois absolu de ${\mathbb Q}^{\hbox{\rm \scriptsize ab}}$ est isomorphe au groupe prolibre de rang $\aleph_0$, $\widehat{F}_\omega$.}, on obtient ce m\^eme r\'esultat en appliquant la proposition \ref{sha} au corps $L={\mathbb Q}^{\hbox{\rm \scriptsize ab}}$. L'int\'er\^et de cette remarque est qu'il n'existe {\it a priori} aucun lien logique entre le $\hbox{\rm PIG}_{/{\mathbb Q}}$ et la conjecture de Shafarevich (voir \cite[Remark 2.2]{DD97b}).
\vskip 2mm 
\noindent
b) $k={\mathbb F}_q(T)$. On applique la proposition \ref{sha} au corps $L=\overline{{\mathbb F}_q}(T)$ : $L$ est la cl\^oture cyclotomique du corps $k$ et est donc inclus dans la cl\^oture r\'esoluble de $k$. Un c\'el\`ebre th\'eor\`eme d\^u \`a Riemann, Harbater et Pop\footnote{Ce th\'eor\`eme affirme que, pour tout corps s\'eparablement clos $K$, le groupe de Galois absolu de $K(T)$ est prolibre de rang \'egal au cardinal de $K$.} assure alors que $L$ a un groupe de Galois absolu prolibre de rang $\aleph_0$. On peut donc en d\'eduire que les groupes profinis qui apparaissent comme groupes de Galois sur ${\mathbb F}_q(T)^{\hbox{\rm \scriptsize r\'es}}$ sont exactement les pro-FnR-groupes de rang $\leq \aleph_0$.
\vskip 2mm 
\noindent
\tab Ce r\'esultat est particuli\`erement int\'eressant quand on le compare au pr\'ec\'edent. En effet, il s'agit de l'analogue classique entre ${\mathbb Q}$ et ${\mathbb F}_q(T)$ : la cl\^oture cyclotomique de ${\mathbb Q}$ est ${\mathbb Q}^{\hbox{\rm \scriptsize ab}}$ et celle de ${\mathbb F}_q(T)$ est $\overline{{\mathbb F}_q}(T)$. Le th\'eor\`eme de Riemann, Harbater et Pop dit alors que l'analogue de la conjecture de Shafarevich est vraie pour ${\mathbb F}_q(T)$ de la m\^eme mani\`ere que le b) dit que l'analogue du a) est vrai pour ${\mathbb F}_q(T)$.
\vskip 10mm 
\noindent
{\large \bf 5.--- Application au Probl\`eme Inverse de Galois Faible.}
\vskip 2mm 
\noindent
\tab Dans cette derni\`ere section, nous donnons tout d'abord des conditions suffisantes sur une classe ${\mathscr C}$ de groupes finis et un corps $k$ pour que le PIGF poss\`ede une r\'eponse positive sur les sous-corps de $k^{\mathscr C}$ contenant $k$. En particulier, dans le th\'eor\`eme \ref{thm clo1} \`a venir, nous g\'en\'eralisons le th\'eor\`eme \ref{thm intro 2} de l'introduction. Nous appliquons ensuite notre \'etude \`a plusieurs classes de groupes finis explicites. Le corollaire \ref{coro P02} qui suit pr\'ecise l'\'ecart existant entre le PIG et sa version Faible annonc\'e dans le th\'eor\`eme \ref{thm intro 3} de l'introduction.
\vskip 5mm 
\noindent
{\bf 5.1.--- R\'esolution du Probl\`eme Inverse de Galois Faible sur certaines cl\^otures d'un corps hilbertien.}
\vskip 2mm 
\noindent
{\bf Th\'eor\`eme\num\label{thm clo1}} {\it On se donne une pr\'e-formation $\mathscr{C}$ et un corps hilbertien $k$ de caract\'eristique $p \geq 0$.
\vskip 2mm
\noindent
{\rm{1)}} Si les conditions
\vskip 1mm
\noindent
\tab $\bullet$ $p \not=2$ ou $\overline{\mathbb{F}_2} \subseteq k$,
\vskip 1mm
\noindent
\tab $\bullet$ il existe une infinit\'e d'entiers $n \geq 1$ tels que le groupe altern\'e $A_n$ n'appartienne pas \`a ${\mathscr C}$,
\vskip 1mm
\noindent
sont satisfaites, alors le $\hbox{\rm PIGF}_{/L}$ admet une r\'eponse positive pour tout corps interm\'ediaire $k\subseteq L\subseteq k^{\tiny \mathscr C}$.
\vskip 2mm 
\noindent
{\rm{2)}} Si les conditions
\vskip 1mm
\noindent
\tab $\bullet$ $p=2$,
\vskip 1mm
\noindent
\tab $\bullet$ il existe une infinit\'e d'entiers impairs $n \geq 1$ tels que $A_n$ n'appartienne pas \`a ${\mathscr C}$,
\vskip 1mm
\noindent
\tab $\bullet$ il existe un groupe fini simple non ab\'elien $G_0$ n'appartenant pas \`a ${\mathscr C}$ et admettant une r\'ealisation r\'e-

guli\`ere sur $k$ ne poss\'edant que $\infty$ comme point de branchement,
\vskip 1mm 
\noindent
sont satisfaites, alors le $\hbox{\rm PIGF}_{/L}$ admet une r\'eponse positive pour tout corps interm\'ediaire $k\subseteq L\subseteq k^{\tiny \mathscr C}$.
\vskip 2mm
\noindent
\tab Plus g\'en\'eralement, dans chacun des cas {\rm 1)} et {\rm 2)} ci-dessus, pour tout groupe fini $G$ et tout corps interm\'ediaire $k \subseteq L \subseteq k^{\mathscr C}$, il existe une suite $(F_n/ L)_{n \geq 1}$ d'extensions finies s\'eparables de groupe d'automorphismes $G$ telle que les corps $F_1, \dots, F_{n}$ soient lin\'eairement disjoints sur $L$ pour tout $n \geq 2$.}
\vskip 2mm 
\noindent
{\bf Preuve :} Notons tout d'abord que, pour qu'il existe, pour tout groupe fini $G$ et tout corps interm\'ediaire $k\subseteq L\subseteq k^{\tiny \mathscr C}$, une suite $(F_n/ L)_{n \geq 1}$ d'extensions finies s\'eparables de groupe d'automorphismes $G$ telle que les corps $F_1, \dots, F_{n}$ soient lin\'eairement disjoints sur $L$ pour tout $n \geq 2$, il suffit de trouver une famille $\{\Gamma_n\}_{n \geq 1}$ de groupes finis non triviaux v\'erifiant les quatre conditions suivantes :
\vskip 2mm 
\noindent
a) $\Gamma_n \in \mathscr{C}^*$ pour tout entier $n \geq 1$,
\vskip 2mm
\noindent
b) pour tout entier $n \geq 1$, aucun quotient simple de $\Gamma_n$ n'est ab\'elien,
\vskip 2mm
\noindent
c) $\{\Gamma_n\}_{n \geq 1}$ v\'erifie la propri\'et\'e de r\'ealisation \hbox{\sc (r\'eal)} (\'enonc\'ee avant le lemme \ref{easy}),
\vskip2mm
\noindent
d) $\Gamma_n$ est groupe de Galois r\'egulier sur $k$ pour tout entier $n \geq 1$.
\vskip 2mm
\noindent
En effet, fixons un groupe fini $G$ et un corps interm\'ediaire $k \subseteq L \subseteq k^\mathscr{C}$. Par c), il existe un entier $n \geq 1$ et un sous-groupe $H$ de $\Gamma_n$ tel que $N_{\Gamma_n}(H)/H \cong G$. Fixons une extension finie galoisienne $E/k(T)$ telle que $E/k$ soit r\'eguli\`ere et v\'erifiant ${\rm{Gal}}(E/k(T))=\Gamma_n$ (cf. d)). Comme $k$ est hilbertien, $\Gamma_n$ est non trivial et $E/k$ est r\'eguli\`ere, il existe une suite $(t_m)_{m \geq 1}$ d'\'el\'ements de $\mathbb{P}^1(k)$ telle que ${\rm{Gal}}(E_{t_m}/k) = \Gamma_n$ pour tout $m \geq 1$ et telle que les corps $E_{t_1}, \dots, E_{t_m}$ soient lin\'eairement disjoints sur $k$ pour tout $m \geq 2$. \'Etant donn\'e $m \geq 1$, le groupe de Galois de l'extension $E_{t_1} \cdot \cdots \cdot E_{t_{m+1}}/k$ s'identifie au produit direct de $m+1$ copies du groupe $\Gamma_n$. Or $\Gamma_n$ est dans ${\mathscr C}^*$  (cf. a)) et ${\mathscr C}^*$ v\'erifie $(\hbox{\rm C}_2)$ (en vertu du 2) de la proposition \ref{transmission 2}). Par cons\'equent, ${\rm{Gal}}(E_{t_1} \cdot \cdots \cdot E_{t_{m+1}}/k)=\Gamma_n^{m+1}$ est dans ${\mathscr C}^{*}$. De plus, par b), aucun quotient simple du groupe $\Gamma_n^{m+1}$ n'est ab\'elien. Par la proposition \ref{propbrubru}, on a donc 
\begin{equation} \label{eq 5.1}
[E_{t_1} \cdot \cdots \cdot E_{t_{m+1}} \cdot L:L]=[E_{t_1} \cdot \cdots \cdot E_{t_{m+1}}:k].
\end{equation}
De m\^eme, on a
\begin{equation} \label{eq 5.2}
[E_{t_1} \cdot \cdots \cdot E_{t_{m}} \cdot L:L]=[E_{t_1} \cdot \cdots \cdot E_{t_{m}}:k]\ \ \hbox{\rm et}\ \ [E_{t_{m+1}} \cdot L:L]=[E_{t_{m+1}}:k].
\end{equation}
Comme les corps $E_{t_1}, \dots, E_{t_{m+1}}$ sont lin\'eairement disjoints sur $k$, on a
\begin{equation} \label{eq 5.3}
[E_{t_1} \cdot \cdots \cdot E_{t_{m+1}}:k]=[E_{t_1} \cdot \cdots \cdot E_{t_{m}}:k][E_{t_{m+1}}:k].
\end{equation}
En combinant \eqref{eq 5.1}, \eqref{eq 5.2} et \eqref{eq 5.3}, on obtient $
[E_{t_1} \cdot \cdots \cdot E_{t_{m+1}} \cdot L:L]=[E_{t_1} \cdot \cdots \cdot E_{t_{m}} \cdot L:L] [ E_{t_{m+1}} \cdot L:L]
$, ce qui montre bien que les corps $E_{t_1} \cdot L, \dots, E_{t_{m+1}} \cdot L$ sont lin\'eairement disjoints sur $L$. Enfin, comme d\'ej\`a not\'e dans \eqref{eq 5.2}, on a ${\rm{Gal}}(E_{t_{m}} \cdot L/L)=\Gamma_n$ pour tout $m \geq 1$. Le lemme \ref{easy} fournit alors la conclusion d\'esir\'ee.
\vskip 2mm
\noindent
\tab Nous d\'emontrons maintenant que, dans chacun des cas 1) et 2), l'on peut toujours trouver une famille $\{\Gamma_n\}_{n \geq 1}$ de groupes finis non triviaux v\'erifiant les quatre conditions ci-dessus. Pour ce faire, nous utiliserons le r\'esultat pr\'eliminaire de th\'eorie des groupes suivant :
\vskip 2mm
\noindent
{\it Si un groupe $\Gamma$ est extension d'un groupe $G_2$ par un groupe $G_1$, alors tout quotient simple de $\Gamma$ est un quotient de $G_2$ ou de $G_1$. En particulier, si un groupe $\Gamma$ est extension d'une puissance d'un groupe simple $G_2$ par un groupe simple $G_1$, alors $G_1$ et $G_2$ sont les seuls quotients simples possibles de $\Gamma$ \footnote{En effet, fixons un groupe simple $S$ et un \'epimorphisme $\pi:\Gamma \longrightarrow S$. Puisque $G_2\trianglelefteq \Gamma$, on a $\pi(G_2)\trianglelefteq S$. Ainsi, soit $\pi(G_2)=S$ et $S$ est alors un quotient de $G_2$, soit $\pi(G_2)=\{1\}$ et, dans ces conditions, on a alors $G_2\leq \hbox{\rm ker}(\pi)$, ce qui fournit un \'epimorphisme $G_1=\Gamma/G_2 \longrightarrow \Gamma/\hbox{\rm ker}(\pi) \cong S$.}.}
\vskip 2mm
\noindent
\tab Supposons tout d'abord que l'on soit dans le cas 1) et $p \not=2$. Fixons un ensemble infini $S$ d'entiers $n \geq 7$ tel que, pour tout $n \in S$, le groupe $A_n$ n'appartienne pas \`a $\mathscr{C}$. Puisque tout groupe fini $G$ se plonge dans $A_n$ pour au moins un entier $n \in S$, l'on voit que, si $m \in S$, la famille des extensions de $A_m^N$ par $A_n$ (index\'ee par $N \geq 1$ et $n \in S$) v\'erifie c) (cf. th\'eor\`eme \ref{thm pigrf}). De plus, la condition d) est v\'erifi\'ee par le lemme \ref{An}. Fixons une telle extension $\Gamma$ de $A_m^N$ par $A_n$. Puisque $A_m$ et $A_{n}$ sont simples et n'appartiennent pas \`a ${\mathscr C}$, on voit que $A_m$ et $A_{n}$ sont \'el\'ements de ${\mathscr C}^*$ et, puisque ${\mathscr C}^*$ v\'erifie $(\hbox{\rm C}_2)$ (en vertu du 2) de la proposition \ref{transmission 2}), on a $\Gamma \in {\mathscr C}^*$. Le r\'esultat pr\'eliminaire ci-dessus assure enfin que tout quotient simple de $\Gamma$ vaut $A_m$ ou $A_n$, qui ne sont visiblement pas ab\'eliens.
\vskip 2mm 
\noindent
\tab Supposons maintenant que l'on soit dans le cas 2). Fixons un groupe fini $G$ et un ensemble infini $S$ d'entiers impairs $n \geq 7$ tel que, pour tout $n \in S$, le groupe $A_n$ n'appartienne pas \`a $\mathscr{C}$. Pour $n_G \in S$ tel que $G \subseteq A_{n_G}$, le groupe $A_{n_G}$ est , en vertu de \cite[Theorem 11]{Bri04} et du lemme \ref{regular}, groupe de Galois r\'egulier sur $k$. Ainsi, en appliquant le cas 1) du th\'eor\`eme \ref{thm pigrf2} et le lemme \ref{regular}, on voit qu'il existe un groupe fini $\Gamma_G$, extension de $G_0^N$ par $A_{n_G}$ pour un certain $N \geq 1$, tel que $\Gamma_G$ soit groupe de Galois r\'egulier sur $k$ et tel que $N_{\Gamma_G}(H)/H \cong G$ pour un certain sous-groupe $H$ de $\Gamma_G$. Ainsi la famille $\{\Gamma_G\}_G$ v\'erifie les conditions c) et d) ci-dessus. On montre alors comme dans le cas pr\'ec\'edent que les conditions a) et b) sont \'egalement v\'erifi\'ees dans cette situation.
\vskip 2mm 
\noindent
\tab Supposons enfin que l'on soit dans le cas 1) et que $\overline{\mathbb{F}_2} \subseteq k$. Alors, d'apr\`es \cite{AOS94}, pour tout $n \geq 7$, le groupe $A_n$ poss\`ede une r\'ealisation r\'eguli\`ere sur $k$ n'admettant que $\infty$ comme point de branchement, et l'on peut donc appliquer le m\^eme raisonnement que celui adopt\'e dans le cas 2).
\fin
\vskip 2mm
\noindent
{\bf Remarques :} 1) Plusieurs groupes finis simples non ab\'eliens poss\`edent une r\'ealisation r\'eguli\`ere sur tout corps de caract\'eristique 2 n'admettant que $\infty$ comme point de branchement. Par exemple, comme d\'ej\`a utilis\'e dans la d\'emonstration du th\'eor\`eme \ref{thm intro 1}, le groupe de Mathieu ${\rm{M}}_{23}$ v\'erifie cette condition. D'autres exemples explicites sont le groupe projectif sp\'ecial lin\'eaire ${\rm{PSL}}_3({\mathbb F}_2)$ (voir \cite[Theorem 5.3]{AY94a}) ou le groupe de Mathieu ${\rm{M}}_{24}$ (cf. \cite{AY94b}). Par contre, il n'est pas clair que, pour tout $n \geq 5$, le groupe altern\'e $A_n$ v\'erifie cette condition\footnote{Il n'est m\^eme pas clair {\it a priori} que tout groupe altern\'e soit groupe de Galois r\'egulier sur tout corps de caract\'eristique 2.}, ceci expliquant l'hypoth\`ese plus restrictive en caract\'eristique 2 dans le th\'eor\`eme \ref{thm clo1} si $k$ ne contient pas ${\overline{{\mathbb F}_2}}$.
\vskip 2mm
\noindent
2) Rappelons qu'un corps $k$ est {\it{RG-hilbertien}} (une terminologie introduite par M. Fried et H. V\"olklein dans \cite{FV92}) si toute extension finie galoisienne $E/k(T)$ telle que $E/k$ soit r\'eguli\`ere poss\`ede au moins une sp\'ecialisation $E_{t_0}/k$ (avec $t_0$ dans $\mathbb{P}^1(k)$) telle que ${\rm{Gal}}(E_{t_0}/k) = {\rm{Gal}}(E/k(T))$. S'il est clair que tout corps hilbertien est n\'ecessairement RG-hilbertien, la r\'eciproque n'est pas vraie en g\'en\'eral, cf. \cite{FV92} et \cite{DH99}. Cependant, la conclusion des cas 1) et 2) du th\'eor\`eme \ref{thm clo1} reste vraie si $k$ est seulement RG-hilbertien\footnote{Par contre, comme il n'est pas clair {\it{a priori}} que toute extension finie \'separable d'un corps RG-hilbertien soit un corps RG-hilbertien (cf. \cite[\S5]{DH99}), nous ne pouvons garantir la conclusion plus g\'en\'erale du th\'eor\`eme \ref{thm clo1} dans cette situation.}.
\vskip 2mm
\noindent
3) En vertu du 2) de la proposition \ref{transmission} (et puisque $\mathscr{C}$ et $\widehat{\mathscr{C}}$ contiennent les m\^emes groupes simples), la conclusion du th\'eor\`eme \ref{thm clo1} est m\^eme vraie pour tout corps interm\'ediaire $k \subseteq L \subseteq k^{\widehat {\mathscr{C}}}$. On voit en fait que la conclusion est vraie pour tout entier naturel $n$ et tout corps interm\'ediaire $k^{\tiny {\mathscr C}[n]} \subseteq L \subseteq (k^{\tiny {\mathscr C}[n]})^{\widehat{\mathscr{C}}}$.
\vskip 2mm
\noindent
4) En utilisant le troisi\`eme point du 1) du th\'eor\`eme \ref{th2brubru} \`a la place de la proposition \ref{propbrubru}, l'on voit que, dans la preuve du th\'eor\`eme \ref{thm clo1}, l'on peut se passer de la condition b) si $\mathscr{C}$ est une pr\'e-vari\'et\'e ou une formation. Comme d\'ej\`a mentionn\'e dans la remarque qui pr\'ec\`ede le th\'eor\`eme \ref{th3brubru}, il en est de m\^eme si $\mathscr{C}$ contient tous les groupes finis d'ordre premier.
\vskip 2mm
\noindent
\tab Comme premi\`ere application du th\'eor\`eme \ref{thm clo1}, on peut consid\'erer la classe ${\mathscr C}(\hbox{\rm r\'es})$ des groupes finis r\'esolubles. C'est une vari\'et\'e extensive qui ne contient aucun groupe fini simple non ab\'elien. Les 1) et 2) du th\'eor\`eme  \ref{thm clo1} fournissent alors l'\'enonc\'e suivant :
\vskip 2mm
\noindent
{\bf Corollaire\num\label{coro res}} {\it Pour tout corps hilbertien $k$, le $\hbox{\rm PIGF}_{/L}$ admet une r\'eponse positive pour tout corps interm\'ediaire $k \subseteq L \subseteq k^{\hbox{\rm \scriptsize r\'es}}$.} 
\vskip 2mm
\noindent
{\bf Remarque :} 1) Comme pour le th\'eor\`eme \ref{1.5 general}, ce r\'esultat est novateur dans la mesure o\`u de nombreux corps interm\'ediaires $k \subseteq L \subseteq k^{\hbox{\rm \scriptsize r\'es}}$ sont non hilbertiens et d\'epassent donc le cadre d'application de \cite{LP18}. Par exemple, dans l'ensemble de ces corps interm\'ediaires non hilbertiens figure, pour tout nombre premier $p$, la pro-$p$-extension maximale de $k$, ce corps \'etant la ${\mathscr C}(\hbox{\rm p})$-cl\^oture de $k$ \footnote{Ce corps est bien non hilbertien en vertu de la remarque qui pr\'ec\`ede la proposition \ref{propbrubru}.}. Comme autres corps interm\'ediaires, l'on trouve \'egalement les cl\^otures ab\'elienne $k^{\rm ab}$ et nilpotente $k^{\rm nil}$ du corps $k$. Ces cas peuvent \^etre obtenus directement par application du th\'eor\`eme \ref{thm clo1} aux classes ${\mathscr C}({\rm ab})$ et ${\mathscr C}({\rm nil})$ \footnote{Par des r\'esultats classiques de Kuyk (cf. \cite[Corollaires 1 \& 2]{Kuy70}), les corps $k^{\rm ab}$ et $k^{\rm nil}$ sont hilbertiens si $k$ l'est. Par cons\'equent, le fait que le PIGF poss\`ede une r\'eponse positive sur ces deux corps, qui est alors une cons\'equence de \cite{LP18}, est r\'eobtenu ici sans utiliser leur caract\`ere hilbertien.}.
\vskip 2mm
\noindent
2) La pro-$2$-extension maximale de $\mathbb Q$ est \'egale au corps $\hbox{\bf C}(i)$, o\`u $\hbox{\bf C}$ d\'esigne le corps des nombres r\'eels constructibles \`a la r\`egle et au compas. Si $L$ d\'esigne un sous-corps de $\hbox{\bf C}$ distinct de $\mathbb{Q}$, alors le corps $L$ n'est la ${\mathscr C}$-cl\^oture de $\mathbb Q$ pour aucune pr\'e-formation ${\mathscr C}$. En effet, si c'\'etait le cas, alors, par le premier point du 1) du th\'eor\`eme \ref{th2brubru}, on aurait ${\mathbb Z}/2{\mathbb Z}\in {\mathscr C}$ et donc $i\in L$, ce qui est impossible. Ceci permet en particulier d'exhiber un exemple d'extension galoisienne $M/{\mathbb Q}$ telle que $M$ ne soit pas hilbertien, telle que le $\hbox{\rm PIGF}_{/M}$ admette une r\'eponse positive et telle que $M$ ne soit la ${\mathscr C}$-cl\^oture de $\mathbb Q$ pour aucune pr\'e-formation ${\mathscr C}$ : la cl\^oture pythagoricienne $\mathbb{Q}^{\rm{pyth}}$ de ${\mathbb Q}$, qui est par d\'efinition le plus petit corps $L$ contenant $\mathbb{Q}$ et tel que $L^2=L^2+L^2$. Par ailleurs, comme annonc\'e dans le r\'esum\'e de cet article, le $\hbox{\rm PIG}_{/\mathbb{Q}^{\tiny \rm pyth}}$ admet une r\'eponse n\'egative puisqu'en vertu du th\'eor\`eme de Diller et Dress le groupe ${\mathbb Z}/4{\mathbb Z}$ n'est pas groupe de Galois sur $\mathbb{Q}^{\rm{pyth}}$ (on renvoie par exemple \`a \cite{Des01} pour un aper\c cu des propri\'et\'es de la cl\^oture pythagoricienne d'un corps donn\'e).
\vskip 2mm
\noindent
\tab Un autre exemple d'application peut \^etre obtenu en consid\'erant, pour un entier $n\geq 1$ donn\'e, la classe ${\mathscr C}(\leq{\rm{n}})$ des groupes finis d'ordre au plus $n$. Il s'agit visiblement d'une pr\'e-vari\'et\'e qui ne contient qu'un nombre fini de groupes altern\'es. De plus, pour $n \leq 244 \, \, 823 \, \, 039 \, \, (=|{\rm{M}}_{24}|-1)$, elle ne contient pas le groupe de Mathieu ${\rm{M}}_{24}$ qui, comme d\'ej\`a rappel\'e, poss\`ede une r\'ealisation r\'eguli\`ere sur tout corps de caract\'eristique $2$ n'admettant que $\infty$ comme point de branchement\footnote{Ce groupe est le plus grand groupe fini simple non ab\'elien v\'erifiant cette propri\'et\'e que nous connaissons.}. Par les 1) et 2) du th\'eor\`eme \ref{thm clo1}, on a donc :
\vskip 2mm
\noindent
{\bf{Corollaire\num\label{coro clo5}}} {\it \'Etant donn\'es un entier $n \geq 1$ et un corps hilbertien $k$ de caract\'eristique $p \geq 0$, si $p \not =2$ ou $\overline{{\mathbb F}_2} \subseteq k$ ou $n \leq 244 \, \, 823 \, \, 039$, alors le $\hbox{\rm PIGF}_{/L}$ admet une r\'eponse positive pour tout corps interm\'ediaire $k \subseteq L \subseteq k^{{\mathscr C}(\leq {\rm n})}$.}
\vskip 5mm
\noindent
{\bf 5.2.--- Sur l'\'ecart entre le Probl\`eme Inverse de Galois et sa version Faible.}
\vskip 2mm
\noindent
\tab Comme d\'ej\`a vu dans le \S4.2, tout groupe fini se r\'ealisant comme groupe de Galois sur la cl\^oture r\'esoluble d'un corps quelconque est n\'ecessairement un groupe fortement non r\'esoluble. Ainsi le corollaire \ref{coro res} fournit un nouvel exemple de corps sur lequel le PIGF admet une r\'eponse positive, mais pas le PIG\footnote{Comme d\'ej\`a mentionn\'e dans l'introduction, le corps $\mathbb{Q}^{\rm{tr}}$ v\'erifie \'egalement cette propri\'et\'e.}. 
\vskip 2mm
\noindent
\tab Un exemple similaire peut \^etre obtenu en consid\'erant, pour un sous-ensemble non vide ${\mathscr P}_0$ donn\'e de l'ensemble ${\mathscr P}$ des nombres premiers, la classe ${\mathscr C}({{\mathscr P}_0})$ des groupes finis $G$ tels que tout diviseur premier de $|G|$ soit dans ${\mathscr P}_0$. Visiblement, ${\mathscr C}({{\mathscr P}_0})$ est une vari\'et\'e extensive non triviale et l'on voit que, si ${\mathscr P}_0\ne {\mathscr P}$, alors ${\mathscr C}({{\mathscr P}_0})$ ne contient qu'un nombre fini de groupes altern\'es. Par cons\'equent, d'apr\`es le 1) du th\'eor\`eme \ref{thm clo1} et le 2) du th\'eor\`eme \ref{th2brubru}, si $k$ est un corps hilbertien de caract\'eristique $p \geq 0$ tel que $p \not=2$ ou $\overline{{\mathbb F}_2} \subseteq k$, alors le $\hbox{\rm PIGF}_{/k^{{\mathscr C}{\rm{(}}{\mathscr P}_0{\rm{)}}}}$ admet une r\'eponse positive, mais aucun \'el\'ement non trivial de ${\mathscr C}({{\mathscr P}_0})$ ne se r\'ealise comme groupe de Galois sur $k^{{\mathscr C}{\rm{(}}{\mathscr P}_0{\rm{)}}}$.
\vskip 2mm
\noindent
\tab Puisque tout groupe fini non trivial $G$ appartient \`a ${\mathscr C}$(${\mathscr P}_0)$ pour au moins un ${\mathscr P}_0$ (d\'ependant de $G$), on obtient le r\'esultat ci-dessous qui g\'en\'eralise le th\'eor\`eme \ref{thm intro 3} \'enonc\'e dans l'introduction :
\vskip 2mm
\noindent
{\bf Corollaire\num\label{coro P02}} {\it On se donne 
\vskip2mm
\noindent
$\bullet$ une famille non vide $\{G_i\}_{i \in I}$ de groupes finis non triviaux telle qu'il existe au moins un nombre premier ne divisant aucun des ordres des groupes $G_i$ ($i \in I$),
\vskip2mm
\noindent
$\bullet$ un sous-ensemble non vide strict ${\mathscr P}_0$  de ${\mathscr P}$ contenant tous les diviseurs premiers des ordres des groupes $G_i$ ($i \in I$),
\vskip2mm
\noindent
$\bullet$ un corps hilbertien $k$ de caract\'eristique $p \geq 0$.
\vskip2mm
\noindent
Supposons $p \not=2$ ou $\overline{{\mathbb F}_2} \subseteq k$. Alors le $\hbox{\rm PIGF}_{/k^{{\mathscr C}{\rm{(}}{\mathscr P}_0{\rm{)}}}}$ admet une r\'eponse positive mais, pour tout $i \in I$, le groupe $G_i$ ne se r\'ealise pas comme groupe de Galois sur $k^{\tiny {\mathscr C}({\mathscr P}_0)}$.}
\vskip 2mm
\noindent
{\bf Remarques :} 1) Comme $I$ est non vide, la remarque qui pr\'ec\`ede la proposition \ref{propbrubru} permet d'affirmer qu'aucun corps $k^{\tiny {\mathscr C}({\mathscr P}_0)}$ comme dans le corollaire \ref{coro P02} n'est hilbertien.
\vskip 2mm
\noindent
2) Comme autre exemple de classes permettant de montrer le th\'eor\`eme \ref{thm intro 3}, on peut citer la famille $\{{\mathscr C}(\hbox{\rm A}_{\geq n})^*\}_{n\geq 5}$, o\`u ${\mathscr C}(\hbox{\rm A}_{\geq n})$ d\'esigne la classe constitu\'ee des groupes altern\'es $A_m$ pour $m\geq n$ et du groupe trivial. En effet, puisque ${\mathscr C}(\hbox{\rm A}_{\geq n})$ est visiblement une pr\'e-formation, on d\'eduit des 1) et 2) de la proposition \ref{transmission 2} que la classe ${\mathscr C}(\hbox{\rm A}_{\geq n})^*$ v\'erifie les conditions $(\hbox{\rm C}_{1,2})$. Par le r\'esultat pr\'eliminaire de th\'eorie des groupes du d\'ebut de la preuve de la proposition \ref{propbrubru}, cette classe est en fait une formation extensive. Si $G$ d\'esigne un groupe fini d'ordre au plus $n$, alors $G\in {\mathscr C}(\hbox{\rm A}_{\geq n})^*$ et donc, en appliquant le 1) du th\'eor\`eme \ref{thm clo1}, le 2) du th\'eor\`eme \ref{th2brubru} et la remarque qui pr\'ec\`ede la proposition \ref{propbrubru}, l'on voit que, si $k$ est un corps hilbertien de caract\'eristique $p \geq 0$ tel que $p \not=2$ ou $\overline{{\mathbb F}_2} \subseteq k$, alors $k^{{\mathscr C}(\hbox{\rm A}_{\geq n})^*}$ est un corps non hilbertien pour lequel le PIGF admet une r\'eponse positive, mais sur lequel $G$ ne se r\'ealise pas groupe de Galois.
\vskip 2mm
\noindent
\tab Cette famille de classes contient les deux familles $\{{\mathscr C}(\leq{\rm{n}})\}_{n \geq 1}$ et $\{{\mathscr C}$(${\mathscr P}_0)\}_{\emptyset \not= \mathscr{P}_0 \subsetneq \mathscr{P}}$ pr\'esent\'ees pr\'ec\'edemment, ainsi que la classe ${\mathscr C}(\hbox{\rm r\'es})$ \footnote{au sens o\`u chacune de ces classes est contenue dans ${\mathscr C}(\hbox{\rm A}_{\geq n})^*$ pour au moins un entier $n \geq 5$.}. Elle est int\'eressante car la suite $\left({\mathscr C}(\hbox{\rm A}_{\geq n})^*\right)_{n\geq 5}$ est croissante pour l'inclusion et, par le raisonnement ci-dessus, on a $\bigcup_{n\geq 5}{\mathscr C}(\hbox{\rm A}_{\geq n})^*={\mathscr C}(\hbox{\rm gr})$. Ainsi, en posant $k_n=k^{{\mathscr C}(\hbox{\rm A}_{\geq n})^*}$ pour tout $n \geq 5$, on obtient une suite croissante de corps $k\subseteq k_5\subseteq k_6\subseteq \cdots \subseteq k_n \subseteq \cdots$ telle que $\bigcup_{n\geq 5}k_n=k^{\hbox{\rm \scriptsize s\'ep}}$ et telle que, pour tout $n\geq 5$, le PIGF admet une r\'eponse positive sur $k_n$, mais pas le PIG.
\vskip 2mm
\noindent
\tab Pour conclure cet article, nous attirons l'attention sur le fait que nous pouvons en fait construire une infinit\'e de r\'ealisations lin\'eairement disjointes d'un groupe fini donn\'e comme groupe d'automorphismes sur la ${\mathscr C}$-cl\^oture consid\'er\'ee du corps hilbertien $k$ dans les corollaires \ref{coro res} \`a \ref{coro P02}, et l'on peut aussi supposer que $k$ est seulement RG-hilbertien.
\bibliography{Biblio2}
\bibliographystyle{alpha}
\vskip 1cm
\noindent
{\bf Bruno Deschamps}
\vskip 2mm
\noindent
{\sc Laboratoire de Math\'ematiques Nicolas Oresme, CNRS UMR 6139}\\
Universit\'e de Caen - Normandie\\
BP 5186, 14032 Caen Cedex - France\\
------------------------------\\
{\sc D\'epartement de Math\'ematiques --- Le Mans Universit\'e}\\
Avenue Olivier Messiaen, 72085 Le Mans cedex 9 - France\\
E-mail : Bruno.Deschamps@univ-lemans.fr
\vskip 5mm
\noindent
{\bf Fran\c cois Legrand}
\vskip 2mm
\noindent
{\sc Institut f\"ur Algebra, Fachrichtung Mathematik}\\
TU Dresden, 01062 Dresden, Germany\\
E-mail : francois.legrand@tu-dresden.de
\end{document}